\documentclass{amsart}
\usepackage{pdfsync}
\usepackage{amsmath}
\usepackage{amssymb}
\usepackage{amsthm}
\usepackage{amscd}
\usepackage{mathrsfs}
\usepackage{xcolor}
\usepackage{ulem}
\usepackage{bm}
\usepackage{graphicx}
\usepackage{longtable}

\newtheorem{thm}{Theorem}
\newtheorem{lem}[thm]{Lemma}

\newtheorem{prop}[thm]{Proposition}
\newtheorem*{PP}{Picard's Problem}
\theoremstyle{definition}
\newtheorem{defin}[thm]{Definition}
\newtheorem{remark}[thm]{Remark}

\newtheorem{example}[thm]{Example}
\newtheorem{conj}[thm]{Conjecture}

\def\bC{{\mathbb C}}

\def\QX{{\mathbb Q}}

\def\NX{{\mathbb N}}

\numberwithin{equation}{section}
\numberwithin{thm}{section}

\newcommand{\Const}{\text{\rm Const}}

\begin{document}
\title[Symbolic Integration of Differential Forms]{Symbolic Integration of Differential Forms:\\ from Abel to Zeilberger}

\thanks{S.\ Chen and Y.\ Wang were partially supported by
National Key R\&D Programs of China (No.~ 2023YFA1009401), the NSFC grants (No.~12271511), the
CAS Funds of the Youth Innovation Promotion Association (No.~Y2022001),
and the Strategic Priority Research Program of the Chinese Academy of
Sciences (No.~XDB0510201).  Y.\ Wang was supported by the Austrian FWF grant 10.55776/I6130.
This work was supported by the International Partnership Program of Chinese Academy of Sciences (No.\ 167GJHZ2023001FN)}.

\author{Shaoshi Chen}
\address{KLMM,\, Academy of Mathematics and Systems Science, Chinese Academy of Sciences, {\it and }School of Mathematics, University of Chinese Academy of Sciences, Beijing, 100190, China}
\email{schen@amss.ac.cn}

\author{David A. Cox}
\address{Department of Mathematics, Amherst College, Amherst, MA 01002, USA}
\email{dacox@amherst.edu}

\author{Yisen Wang}
\address{KLMM,\, Academy of Mathematics and Systems Science, Chinese Academy of Sciences, {\it and} School of Mathematics, University of Chinese Academy of Sciences, Beijing, 100190, China, {\it and} Johann Radon Institute for Computational and Applied Mathematics,\, Austrian Academy of Sciences, Linz, 4040, Austria} 
\email{wangyisen@amss.ac.cn}

\begin{abstract}
This paper focuses on  symbolic integration of differential forms, with a particular emphasis on historical and modern developments, from Abel's addition theorems for Abelian integrals to Zeilberger's creative telescoping for parameterized integrals. It explores closed rational $p$-forms and provides algorithmic approaches for their integration, extending classical results like Hermite reduction and Liouville's theorem. The integration of closed differential forms with parameters is further examined through telescopers, offering a unified framework for handling both algebraic and transcendental cases.
\end{abstract}


\maketitle

\section{Introduction}
\label{SECT:intro}
Differential forms play a central role in mathematics in both theoretical and computational aspects. They provide an elegant language to unify different kinds of integrals in multivariable calculus and are used in differential geometry~\cite{Spivak},
algebraic topology~\cite{BottTu}, and algebraic geometry~\cite{Griffiths1978} whenever integrals are involved. Symbolic integration of univariate functions has been well-developed starting with the celebrated work by Risch~\cite{Risch1969} and summarized comprehensively in Bronstein's book~\cite{Bronstein} and the {foundational references collected in}~\cite{RaabSinger2022}. A long-standing challenging project in symbolic computation is developing symbolic algorithms for computing integrals of multivariate functions.  In the multivariate setting, differential forms provide a natural language for studying integration problems. 
Symbolic computation of differential forms has been studied in computational algebraic geometry~\cite{Burgisser2007, Burgisser2009}
with applications to computing periods in~\cite{Lairez2016}.
Packages in computer algebra systems such as Maple and Mathematica rely on a variety of  heuristic techniques to integrate differential forms. The goal of this paper is to enrich the algorithmic methods available for symbolic integration of differential forms, 
motivated by Abel's classical addition theorem on integrals of algebraic functions. 

In 1826, Abel submitted a remarkable memoir~\cite{Abel1} to the Paris Academy, who did not publish it until 1841, long after his death. His main focus is on what we now call \textit{Abelian integrals}, which are integrals of the form
\[
\int \!f(x,y)\,dx,
\]
where $x$ and $y$ are related by a polynomial equation $\chi(x,y) = 0$.   Abel considers an auxiliary equation $\theta(x,y) = 0$ that depends on parameters $a, a', a'', \dots$, so that the solutions of the simultaneous equations
\[
\chi(x,y) = \theta(x,y) = 0
\]
are algebraic functions of the parameters.  Abel focuses on the nonconstant solutions and writes them as  $(x_i,y_i)$, $i \in \{1,\dots,\mu\}$.  In Section 2 of~\cite{Abel1}, Abel states the following theorem:
\begin{quote}
	I say that if we denote by $f(x,y)$ some rational function of $x$ and $y$, and if we set
	\[
	dv = f(x_1,y_1)\,dx_1 +  f(x_2,y_2)\,dx_2 + \cdots +  f(x_\mu,y_\mu)\,dx_\mu,
	\]
	the differential $dv$ will be a rational function of the quantities $a, a', a''$, etc.
\end{quote}
This is one of several results in~\cite{Abel1} that go under the name ``Abel's Theorem'' (see~\cite{Kleiman}).  Abel's first proof of the rationality of $dv$  is very terse.  A more detailed argument is given in his last paper~\cite{Abel2}, and modern explanations can be found in~\cite{Cox2024},~\cite[Ch.\ 9]{Edwards}, and~\cite{Griffiths1976}.  In Section 4 of~\cite{Abel1}, Abel notes that this method is ``in general very long, and for functions that are a little complicated, almost impractical.'' He then gives a direct proof of the rationality of $dv$ by expressing it as an explicit rational differential of the parameters.  His astonishing formula is described in~\cite{Cox2024}.

Later in Section 2 of~\cite{Abel1}, Abel notes that:
\begin{quote}
	If now $dv$ is a rational function of the quantities $a, a', a''\dots$ its integral or the quantity $v$ will be an algebraic and logarithmic function of $a, a', a''\dots$\footnote{While Abel says ``algebraic and logarithmic'' in his result, Theorem~\ref{THM:1form} shows that ``algebraic'' can be replaced with ``rational'' (this is also proved in~\cite{Cox2024}). Abel's preference for ``algebraic'' is related to a change made in Section 6 of~\cite{Abel1}.  If $\alpha$ is the number of parameters $a,a',a'',\dots$ and $(x_i,y_i)$, $i \in\{1,\dots,\mu\}$, are the solutions considered above, then Abel claims that $a,a',a'',\dots$ are algebraic functions of the first $\alpha$ of $x_1,\dots,x_\mu$.  This makes $v$ an algebraic and logarithmic function of $x_1,\dots,x_\alpha$.}
\end{quote}
For us, this is where things get interesting. In going from $dv$ to $v$, Abel is clearly referring to the method of partial fractions, which is where the logarithms come from.  But this method applies to rational $1$-forms in one variable, while here, $dv$ is a rational $1$-form in the parameters $a, a', a'',\dots$.  What's missing is the recognition that $dv$ is \textit{closed} (which is easy to prove) and that one needs a result like Theorem~\ref{THM:1form} to guarantee that $v$ is a rational and logarithmic function of the parameters.  Abel avoids this difficulty since his explicit formula for $dv$ leads immediately to an explicit formula for $v$.

In contrast, we show in Theorem~\ref{THM:1form} that this is the case for \textit{any} closed rational $1$-form.  A natural question is whether Theorem~\ref{THM:1form} still holds for closed rational $p$-forms. This is answered affirmatively by Theorem~\ref{THM:pform} in Section~\ref{SECT:pform}.
Deligne's logarithmic comparison theorem~\cite[Prop.~3.1.8]{Deligne1971} shows that, when the polar divisor has normal crossings, differential forms with logarithmic poles suffice to compute the cohomology while preserving the polar locus. By contrast, Theorem~\ref{THM:pform} does not require the polar divisor to have normal crossings, but it may introduce new denominators.

Liouville's  classical theorem is crucial for determining whether an elementary function has an elementary integral. The first step towards symbolic integration of differential forms with elementary-function coefficients would be establishing Liouville's theorem in this setting. In this direction, the first work was done by Caviness and Rothstein~\cite{Caviness1975} who proved a version of Liouville's theorem for integration
in finite terms of line integrals. Using the techniques introduced in Section~\ref{SECT:1form}, we provide a more direct proof in the language of differential $1$-forms. It is still challenging to extend this theorem to the case of $p$-forms.

Symbolic computation of parameterized integrals is an active topic~\cite{Raab2016} with rich applications to Gauss--Manin connections and Picard--Fuchs equations in algebraic geometry~\cite{Manin1958,Lairez2016} and Feymann diagrams in mathematical physics~\cite{muller-weinzierl-zayadeh, Schneider2021}. A powerful tool for evaluating parameterized integrals is Zeilberger's method of creative telescoping~\cite{AlmkvistZeilberger1990, Zeilberger1991}.  The formulation of creative telescoping can be generalized from functions to differential forms.  Given a differential form $\omega$ in $x_1, \ldots, x_m$ with a parameter $t$, a linear differential operator $L$ with respect to $t$ is called a \textit{telescoper} for $\omega$ if $L(\omega)$ is exact over the ground field. The existence problem of telescopers was first studied by Manin in~\cite{Manin1958} for the case of $m$-forms with algebraic functions as coefficients and recently was investigated in~\cite{ChenFengLiSinger2014, ChenFengLiSingerWatt2024} for the case of differential forms with D-finite functions as coefficients. We will present an algorithm {that uses Hermite reduction to compute minimal telescopers of rational differential 1-forms with one parameter.}

The outline of the paper is as follows. After establishing some notation and preliminary lemmas in Section~\ref{SECT:prelim}, we generalize the 1-forms $dv$ studied by Abel to the case of closed rational $1$-forms and develop the corresponding Hermite reduction in Section~\ref{SECT:1form}. The case of closed rational $p$-forms is considered in Section~\ref{SECT:pform}.  In Section~\ref{SECT:Picard}, we recall Picard's problem about  
the exactness of differential $3$-forms over the rational function field and Griffiths--Dwork reduction for solving the smooth case together with an example and a conjecture related to Fermat hypersurfaces.  In  Section~\ref{SECT:Liouville}, we generalize Liouville's theorem to differential $1$-forms. We conclude the paper by studying Zeilberger's creative telescoping for closed rational $1$-forms in Section~\ref{SECT:ZeilCT}.

\section{Preliminaries}
\label{SECT:prelim}
In this section, we first recall some basic facts about differential fields and their extensions~\cite{Bronstein}.  We also discuss differential forms.

\begin{defin}
Let $R$ be a ring (resp.\ field). A \textit{derivation} on $R$ is a map $D:R  \rightarrow R$ such that for any $a,b \in R$, we have $D(a+b)= D(a)+D(b)$ and $D(ab)=aD(b)+D(a)b$. The pair $(R,D)$ is called a \textit{differential ring} (resp.\ \textit{differential field}). The set $\Const_D(R) := \{r\in R\mid D(r)=0\}$ is a subring (resp.\ subfield) of $R$. A pair $(R^*, D^*)$ is called a \textit{differential extension} of $(R,D)$ if $R\subseteq R^*$ and $D^*{\mid_{R}} = D$.
\end{defin}

The following fact is used frequently in symbolic integration~\cite{Bronstein}.

\begin{lem}
\label{LEM:fact}
For a field $F$ of characteristic zero, let $F(y)$ be the field of rational functions in $y$ over $F$. Let $D_y$ denote the derivative $d/dy$ on $F(y)$ satisfying $D_y(y) = 1$ and $D_y(c)=0$ for all $c\in F$.  Given pairwise coprime polynomials $u_1, \ldots, u_n\in F[y]\setminus F$, constants $c_1, \ldots, c_n\in F$, and a rational function $v\in F(y)$, suppose that
\[
\sum_{i=1}^n c_i \frac{D_y(u_i)}{u_i} + D_y(v) = 0.
\]
Then $c_1 = \cdots = c_n = 0$ and $v\in F$.
\end{lem}

\begin{proof}
Fix $i$ and let $\beta$ be a root of $u_i$ of multiplicity $\ell$ in some differential extension of $F$.   By hypothesis, $u_j(\beta) \ne 0$ for $j \ne i$.  A standard calculation shows that
\[ \text{residue}_{y=\beta}\Big( \frac{D_y(u_i)}{u_i}\Big) = \ell \ne 0.\]
The residue of the derivative of a rational function always vanishes, so taking the residue at $y=\beta$ of each side of the equation in the lemma implies $c_i \ell = 0$, and $c_1 = \cdots = c_n = 0$ follows.  Then $D_y(v) = 0$, which implies $v \in F$.
\end{proof}

For a field $k$ of characteristic zero, let $K:=k(x_1, \ldots, x_m)$ be the field of rational functions in $x_1, \ldots, x_m$ over $k$.  Let $\partial_i$ denote the usual partial derivative $\partial/\partial x_i$ that satisfies $\partial_i \partial_j(f) = \partial_j\partial_i(f)$ for all $f\in K$ and $\partial_{i}(c)=0$ for all $c \in k(x_1, \ldots, x_{i-1}, x_{i+1}, \ldots, x_m)$.
Let $\mathcal{U}$ be the universal differential extension of $K$ in which any system of
algebraic differential equations has a solution in $\mathcal{U}$ if it has a solution in some differential extension of $K$ (see~\cite[Ch.~III, Section 7]{kolchin} for the existence of such universal fields). Note that the $\partial_i$'s also commute on $\mathcal{U}$.

We now recall some basic notions about exterior algebras over $\mathcal{U}$ from Lang's book~\cite[Ch.~XIX]{SergeAlgebra}. Given a vector space ${\mathcal{M}}$ 
over $\mathcal{U}$ of dimension $m$ with basis $\{\mathbf{a}_1,\ldots,\mathbf{a}_m\}$, let
\[
	\bigwedge \mathcal{M} = \bigoplus_{p=0}^m {\bigwedge}^p {\mathcal{M}}
\]
be the exterior algebra of $\mathcal{M}$, where $\bigwedge^p \mathcal{M}$ denotes the vector space over $\mathcal{U}$ with basis $\{\mathbf{a}_{i_1} \wedge \cdots \wedge \mathbf{a}_{i_p} \ | \ 1 \leq i_1 < \cdots < i_p \leq m\}$. We call an element in $\bigwedge^p \mathcal{M}$ a $p$-form, and an element $f \in \mathcal{U}$ is called a  $0$-form. Let $d \colon \mathcal{U} \rightarrow \mathcal{M}$ be the map defined by
\[
	d\hskip.01pt f = \partial_1 (f) \hskip1pt \mathbf{a}_1 + \cdots + \partial_m(f)\hskip1pt  \mathbf{a}_m.
\]
Hence $dx_i = \mathbf{a}_i$ by taking $f= x_i$. For the remainder of this paper, we will use $\{ dx_1, \ldots, dx_m\}$ instead of $\{\mathbf{a}_1,\ldots, \mathbf{a}_m\}$. We call $\bigwedge \mathcal{M}$ the \textit{space of differential forms}.
Suppose that $\omega \in {\bigwedge}^p \mathcal{M} $ is written as  
\begin{equation}\label{EQ:pform}
    \sum_{1 \leq i_1 < \cdots < i_p \leq m} f_{i_1,\ldots,i_p} \,dx_{i_1}\wedge \cdots \wedge dx_{i_p}. 
\end{equation}
We call each $f_{i_1,\ldots,i_p}$ a coefficient of~$\omega$.

We recall some properties of operations on differential forms from~\cite{Weintraub,Cartan}. The wedge product $\omega  \wedge \eta$ is bilinear, anticommutative, and associative as follows:
\begin{align*}
	f_{I}\, dx_{I} \wedge g_{J} \, dx_{J} & = f_{I} g_{J} \, dx_{I}\wedge dx_{J}\\
	(f_{I}\,dx_{I}+f'_{I'}\,dx_{I'}) \wedge g_{J}\, dx_{J} &= f_{I}g_{J} \, dx_{I}dx_{J} + f'_{I'}g_{J} \, dx_{I'} dx_{J}\\
	f_{I} \, dx_{I} \wedge g_{J}\, dx_{J} & = (-1)^{|I|\cdot|J|} g_{J} \, dx_{J} \wedge f_{I} \, dx_{I}\\
	(f_{I} \, dx_{I} \wedge g_{J}\, dx_{J}) \wedge h_{R}\, dx_{R} & = f_{I} \, dx_{I} \wedge ( g_{J}\, dx_{J} \wedge h_{R}\, dx_{R}).
\end{align*}
Here,
\[
dx_I:=dx_{i_1}\wedge\cdots\wedge dx_{i_p},
\]
where $I=(i_1,\ldots,i_p)$ is an index tuple. We often abbreviate $dx_I$ as $dx_{i_1}\cdots dx_{i_p}$, with the wedge products understood. Moreover, $|I|=p$ denotes the length of $I$.
If $dx_I$ is any string with some $dx_i$ appearing more than once, then $dx_I=0$.
If $i_1,\ldots,i_k$ are distinct and $\sigma \colon \{i_1,\ldots,i_k\} \rightarrow \{i_1,\ldots,i_k\}$ is a permutation, then
\[dx_{\sigma(i_1)}\cdots dx_{\sigma(i_k)}=\mathrm{sign}(\sigma)\,dx_{i_1}\cdots dx_{i_k},\]
where $\mathrm{sign}(\sigma)=\pm 1$ is the sign of the permutation.
The map $d$ can be extended to a graded derivation on $\bigwedge \mathcal{M}$. 
If $\omega \in \bigwedge^p \mathcal{M}$ is written in the form~\eqref{EQ:pform}, then its \textit{exterior derivative} $d\omega$ is defined as 
\[
     \sum_{1 \leq i_1 < \cdots < i_p \leq m} 
     d\bigl(f_{i_1,\ldots,i_p}\bigr) \,dx_{i_1} \cdots dx_{i_p} .
\]
%
The difference between differential forms in $\mathcal{M}$ and K\"ahler differentials is well explained in~\cite{Li2011}.
\begin{remark}\label{RMK:pullback}
Let $\phi \colon K \rightarrow K$ be an automorphism that is the identity on~$k$. Then $\phi$ extends to an automorphism of~$\mathcal{U}$, which in turn gives a map~$\phi^* \colon \mathcal{M} \rightarrow \mathcal{M}$ defined by~$\phi^*(f \, dx_i) = \phi(f) \, d(\phi(x_i))$. This extends to an automorphism~$\phi^* \colon \bigwedge \mathcal{M} \rightarrow \bigwedge \mathcal{M}$. We say that $\phi^*(\omega)$ is the \textit{pullback} of $\omega$ via~$\phi$. Note that pullbacks commute with the wedge product~$\wedge$ and commute with the exterior derivative~$d$ via the chain rule.
\end{remark}
\begin{defin}
	Let~$\omega \in \bigwedge \mathcal{M} $ and $R$ be a subring of~$\mathcal{U}$. We say that $\omega$ is \textit{closed} if~$d \omega = 0 $ and $\omega $ is \textit{exact} over~$R$ if there exists~$\eta \in \bigwedge \mathcal{M}$ with coefficients in~$R$ such that~$d\eta = \omega$. We call~$\eta$ a primitive of~$\omega$.
\end{defin}

Since $d(d\eta) = 0$, every exact differential form is closed.

We now define operators $d_s$ and $d^s$ that will be used in our recursive arguments.  For any $f\in \mathcal{U}$ and $s\in \{1, 2, \ldots, m\}$, define
\[
d^s(f \, dx_I) = \partial_s(f)\,dx_s\wedge dx_I = \partial_s(f)\,dx_s dx_I.
\]
The definition of $d^s$ extends linearly to all of $\bigwedge \mathcal{M}$.  Then define
\[
d_s = d^1 + \cdots + d^s,
\]
and note that
\[
d = d^1 +\cdots + d^m = d_{m-1} + d^m.
\]

Given $\omega \in \bigwedge^p \mathcal{M}$, we can always decompose $\omega$ into
\[
\omega = \omega_{m-1} + \omega^m,
\]
where all summands $f_{i_1, \ldots, i_p} \, dx_{i_1}\cdots dx_{i_p}$ of $\omega_{m-1}$ do not involve $dx_m$ and $\omega^m = \mu\, dx_m$ with $\mu \in \bigwedge^{p-1} \mathcal{M}$ and free of $dx_m$. The following lemma is an easy consequence of the above definitions.
\begin{lem}
\label{LEM:form}
If $\omega = \omega_{m-1} + \omega^m  \in \bigwedge^p \mathcal{M}$ is as above, then we have
\begin{itemize}
\item[{\rm(i)}] $d^m (\omega^m) = 0$.
\item[{\rm(ii)}] $d(\omega \,dx_m) = d_{m-1}(\omega)\, dx_m$.
\item[{\rm(iii)}] $d\hskip1pt\omega =0$ if and only if  $d_{m-1}(\omega_{m-1}) = 0$ and $d^m(\omega_{m-1})  + d_{m-1}(\omega^m)=0$.
\item[{\rm(iv)}] If $\omega$ is free of $dx_m$, then $d\hskip1pt\omega = 0$ if and only if  $d_{m-1}(\omega) = 0$ and $d^m (\omega) = 0$.
\item[{\rm(v)}] If $\omega $ is free of $x_m$ and $dx_m$, then $\omega = d\mu$ for some $(p-1)$-form $\mu$ if and only if $ \omega = d_{m-1}(\tilde \mu)$ for some $(p-1)$-form $\tilde \mu$ also free of $x_m$ and $dx_m$.
\end{itemize}
\end{lem}

\section{Integration of Closed Rational $1$-Forms}
\label{SECT:1form}

Let $k$ be a field of characteristic zero. 
We first study the integration of closed $1$-forms over~$K=k(x_1,\ldots,x_m)$ allowing logarithmic terms in their antiderivatives.
As formulated by Rosenlicht~\cite{Rosenlicht1975},
a \textit{logarithm} of $f \neq 0$, denoted by~$\log f$, in a $\{ \partial_1,\ldots, \partial_m\}$-field extension of $K$ is an element satisfying~$\partial_i(\log f) = \partial_i(f)/f$ for each~$i \in \{1,\ldots,m\}$.

 Let $F := k(x_1, \ldots, x_{m-1})$ and $\overline{F}$ denote the algebraic closure of~$F$. Assume that $\omega = f_1\, dx_1 + \cdots + f_m\, dx_m$ is a closed rational $1$-form with $f_1,\dots,f_m \in K$.
By the classical integration formula for rational functions in the variable $x_m$ (i.e., partial fractions with respect to $x_m$), we have
\begin{equation}\label{EQ:fm}
f_m = \partial_{m}(g_m)\quad  \text{with}\quad  g_m = a_{m} + \sum_j c_{m, j} \log b_{m, j},
\end{equation}
where  $a_{m} \in K$, $c_{m, j} \in \overline{F}$ and $b_{m, j} \in F(c_{m, j})[x_m]$ 
has positive degree in $x_m$.
Since $\omega$ is closed, we have $\partial_i(f_m) = \partial_m(f_i)$ for~$i\in\{1, \ldots, m-1\}$. Thus
\begin{equation}\label{EQ:partialfm} 
\begin{aligned}
 \partial_m(f_i) &= \partial_i(f_m) = \partial_i(\partial_m(g_m))= \partial_m(\partial_i(g_m))\\
  &= \partial_m\bigg(\partial_i(a_m) + \sum_j \Big(c_{m, j}\frac{\partial_i(b_{m, j})}{b_{m, j}} + \partial_i(c_{m, j}) \log b_{m, j}\Big)\!\bigg) \\
 &= \partial_m\bigg(\partial_i(a_m) + \sum_j c_{m, j} \frac{\partial_i(b_{m, j})}{b_{m, j}} \bigg) + \sum_j \partial_i(c_{m, j})\frac{\partial_m(b_{m, j})}{b_{m, j}},
\end{aligned}
\end{equation}
where the last line follows because $c_{m, j} \in \overline{F}$ implies that $\partial_m(\partial_i(c_{m, j})) = 0$.

For fixed $j$, Lemma~\ref{LEM:fact} implies that $\partial_i(c_{m, j}) = 0$ for all $i\in\{1, \ldots, m-1\}$. This means that each $c_{m, j}$ is a constant in $\overline{F}$, so that $c_{m,j} \in \bar k$. Let $\widetilde{\omega} = \omega - dg_m$. Then $\widetilde{\omega} = \tilde{f}_1 \,dx_1 + \cdots + \tilde{f}_{m-1}\,dx_{m-1}$, where for each $i\in\{1, \ldots, m-1\}$,
\begin{equation*}
\tilde{f}_i = f_i - \partial_i(g_m) = f_i -\partial_i(a_m) - \sum_j c_{m, j} \frac{\partial_i(b_{m, j})}{b_{m, j}}.
\end{equation*}

We next observe that the sum $\sum_j c_{m, j} \frac{\partial_i(b_{m, j})}{b_{m, j}}$  lies in $K$.  As written, it lies in the extension $K_1$ of $K$ obtained by adjoining the $c_{m,j}$.  But the $c_{m,j}$ are the roots of a Rothstein--Trager resultant, and the gcd construction of $b_{m,j}$ in the Rothstein--Trager Algorithm (see Section 2.4 of \cite{Bronstein}) shows that any permutation of the $c_{m,j}$ also permutes the $b_{m,j}$.  It follows that the sum $\sum_j c_{m, j} \frac{\partial_i(b_{m, j})}{b_{m, j}} \in K_1$ is invariant under all permutations of the $c_{m,j}$ and hence lies in the subfield $K$.  Thus $\widetilde{\omega} = \tilde{f}_1 \,dx_1 + \cdots + \tilde{f}_{m-1}\,dx_{m-1}$, where $\tilde{f}_i \in K$ for $i\in\{1, \ldots, m-1\}$.

Since $\omega$ is closed, so is $\widetilde{\omega} = \omega - dg_m$, which implies that
\[
0 = d\hskip1pt\widetilde{\omega} = \sum_{i=1}^{m-1} d\tilde{f}_i \wedge dx_i = \sum_{i=1}^{m-1} \sum_{j=1}^{m-1} \partial_j(\tilde{f}_i) \hskip1pt \, dx_j dx_i + \sum_{i=1}^{m-1} \partial_m(\tilde{f}_i) \hskip1pt \, dx_mdx_i .
\]
Thus $\partial_m(\tilde f_i)=0$ for $i\in\{1, \ldots, m-1\}$. So $\widetilde{\omega}$ is a differential form that is free of $dx_m$ and has coefficients in $F = k(x_1,\dots,x_{m-1})$.  Repeating this process recursively, we obtain the following result.

\begin{thm}
\label{THM:1form}
Let $\omega = f_1\, dx_1 + \cdots + f_m \, dx_m$ be a closed $1$-form with $f_i \in K = k(x_1,\dots,x_m)$. 
Then $\omega = dg$ for some $g$ of the form
\[
g = a + \sum_{i=1}^m \sum_j c_{i,j} \log b_{i,j},
\]
where  $a \in K$, $c_{i,j} \in  \bar k$, and $b_{i,j} \in k(x_1,\ldots,x_{i-1},c_{i,j})[x_i]$ has positive degree in~$x_i$.
\end{thm}

\begin{remark}
Note also that a version of Theorem~\ref{THM:1form} appears as Theorem~8 in~\cite{ChenKoutschan}. In~\cite[(1.5)]{Griffiths1976}, Griffiths states a version of Theorem~\ref{THM:1form} and provides a proof in the case of two variables.
\end{remark}

Ostrogradsky~\cite{Ostrogradsky1845} and Hermite~\cite{Hermite1872} introduced 
a reduction method, called \textit{Hermite reduction} in Bronstein's book~\cite{Bronstein}, that decomposes $f\in k(x)$ uniquely as
\[f = \frac{dg}{dx}+ \frac{a}{b},\]
where $g \in k(x)$ and $a, b\in k[x]$ are such that $\gcd(a, b)=1$, $\deg(a)<\deg(b)$ and $b$ is squarefree. Moreover, $f$ is integrable in $k(x)$, i.e., $f= \frac{dh}{dx}$ for some $h\in k(x)$ if and only if~$a/b=0$. Hermite reduction separates the rational part of the integral of a rational function from  the logarithmic part that is transcendental over $k(x)$. Readers can refer to~\cite{Hermite1872} and Section~2.2 in~\cite{Bronstein} for more details. This reduction has been extended to more general classes of functions, including algebraic functions~\cite{Trager1984,ChenKauersKoutschan2016}, hyperexponential functions~\cite{bostan-chen-chyzak-li-xin}, and D-finite functions~\cite{bostan-chyzak-lairez-salvy, vdH2021, Chen2023}.  

We will extend Hermite Reduction in a different direction, namely to closed rational $1$-forms.  We begin with some preliminaries.


    
\begin{defin}
A rational function $f \in {\bar k}(x_1,\ldots,x_m)$ is said to be \textit{proper in~$x_i$} if $f$ can be written as $a/b$ with $a,b \in {\bar k}[x_1,\ldots,x_m]$ such that $\gcd(a,b)=1$ and $\deg_{x_i}(a) < \deg_{x_i}(b)$. 
\end{defin} 

\begin{remark}\label{RMK:proper} 
By convention, we define~$\deg_{x_i}(0) := - \infty$ for each~$i\in\{1,\ldots,m\}$. Hence $0$ is proper in~$x_i$.  
Also, note that for any nonzero $f \in {\bar k}(x_1,\ldots,x_m)$, $\partial_i(f)/f$ is proper in $x_i$ and has a squarefree denominator in ${\bar k}[x_1,\dots,x_m]$.
\end{remark}
\begin{defin}
Let $\omega$ be a $1$-form written as 
\[
\omega = f_1 \,dx_1+\cdots+f_m \,dx_m,
 \text{ where }f_i\in K.
\]
We call $\omega$ a \textit{simple form} if for each $i \in \{1, \ldots, m\}$, the coefficient $f_i$ is proper in~$x_i$ and its denominator is squarefree in $k[x_1,\ldots,x_m]$.
\end{defin}

\begin{remark}\label{RMK:sumsimple}  Note that the sum of two simple forms is again a simple form since proper rational functions are closed under addition and the least common multiple of squarefree polynomials in $k[x_1,\ldots,x_m]$ is still squarefree. 
\end{remark}

\begin{prop}\label{PROP:simpleform}
Let $\omega$ be a simple form over~$K$. Then $\omega$ is exact over $K$ if and only if~$\omega =0$.
\end{prop}
\begin{proof}
    The sufficiency follows from the equality~$\omega = 0= d(1)$. For the necessity, assume that~$\omega=f_1 \, dx_1 + \cdots + f_m \, dx_m = dg$ for some~$g \in K$. Hence $f_i = \partial_i(g)$ for each~$i$. Since $f_i$ is proper in $x_i$ and has a squarefree denominator, the uniqueness property of Hermite reduction implies that $f_i = 0$. Hence~$\omega = 0$.
\end{proof}

We now present an algorithm for decomposing a closed $1$-form into the sum of an exact $1$-form and a simple form that leaves the base field unchanged. This extends Hermite reduction from univariate rational functions to closed rational  $1$-forms. 

\medskip
\noindent
\textbf{Algorithm} \texttt{HermiteOneForm:}
	Given a closed $1$-form $\omega = f_1\, dx_1 + \cdots + f_m\, dx_m$ with $f_i \in K = k(x_1,\dots,x_m)$, return $g \in K$ and a simple form $r$ with coefficients in $K$ such that
\[\omega = dg + r.\]
\begin{enumerate}
\item When $m=1$, apply Hermite reduction to $f_1$ to obtain $f_1 = \partial_1(g) + h$. Return $g$ and~$r := h\, dx_1$.
\item When $m>1$, apply Hermite reduction to $f_m$ with respect to $x_m$ to obtain
\[f_m = \partial_m(a_m) + \bar{f}_m, \]
where $a_m, \bar{f}_m \in K$, and $\bar{f}_m$ is proper in $x_m$ with a denominator that is squarefree in $F[x_m]$, $F =k(x_1,\dots,x_{m-1})$.
\item Apply the Rothstein--Trager Algorithm to $\bar f_m$ to construct  $c_{m,j} \in \bar k$ and $b_{m, j} \in F(c_{m,j})[x_m]$ of positive degree in $x_m$ such that 
\[\bar{f}_m = \sum_j c_{m, j} \frac{\partial_m(b_{m, j})}{b_{m,j}}.\]
\item For each $i\in\{1,\dots,m-1\}$, compute  $\bar f_i := \sum_j c_{m, j} \frac{\partial_i(b_{m, j})}{b_{m,j}}$ as an element of $K$ and set  $\bar r_m := \bar f_1\,dx_1 + \cdots + \bar f_{m}\,dx_{m}$ and $\omega_{m-1} = \omega - da_m - \bar r_m$. 
\item Apply \texttt{HermiteOneForm} to ${\omega}_{m-1}$ with $m-1$ in place of $m$ to obtain
\[{\omega}_{m-1}  = dg_{m-1}+ r_{m-1},\]
where $g_{m-1}\in F$ and $r_{m-1}$ is a simple form over $F$.
\item 
Return $g := a_m+g_{m-1}$ and $r := \bar{r}_m + r_{m-1}$.
\end{enumerate}

\begin{thm}
Algorithm {\tt HermiteOneForm} is correct.
\end{thm}

\begin{proof} We focus on the case $m > 1$.  In Step (2), we write $f_m = \partial_m(a_m) + \bar f_m$ with $a_m, \bar{f}_m \in K$ and $\bar{f}_m$  proper in $x_m$, and in Step (3), $\bar{f}_m = \sum_j c_{m, j} \frac{\partial_m(b_{m, j})}{b_{m,j}} = \partial_m\big(\sum_j c_{m, j} \log b_{m, j}\big)$.
The discussion following \eqref{EQ:partialfm} shows that the $c_{m,j}$ lie in $\bar k$.  

In Step (4), we define $\bar f_i := \sum_j c_{m, j} \frac{\partial_i(b_{m, j})}{b_{m,j}} = \partial_i\big(\sum_j c_{m, j} \log b_{m, j}\big)$ for all $i \in \{1,\ldots,m-1\}$.   Note that $\bar f_i \in K$ by the proof of Theorem~\ref{THM:1form}. Furthermore,
for any~$1 \leq i < \ell \leq m$,
\begin{equation}
\partial_\ell(\bar f_i) = \partial_\ell\partial_i\big({\textstyle\sum_j} c_{m, j} \log b_{m, j}\big) = \partial_i\partial_\ell\big({\textstyle\sum_j} c_{m, j} \log b_{m, j} \big) = \partial_i(\bar f_\ell).
\end{equation}
So the $1$-form $\bar{r}_m := \bar{f}_1 \, dx_1 + \cdots + \bar{f}_m \, dx_m$ defined in Step (4) is closed. 
Remark~\ref{RMK:proper} also implies $\bar f_i$ is proper in $x_i$ with squarefree denominator for $i \in \{1,\dots,m\}$.  It follows that $\bar{r}_m$ is a simple form with coefficients in $K$.


The form $\omega_{m-1}:=\omega-da_m-\bar{r}_m$ in Step~(4) is closed since~$d a_m$ is exact and both $\omega$ and $\bar{r}_m$ are closed. (Note that $\omega_{m-1}$ is the form $\widetilde \omega$ used in the proof of Theorem~\ref{THM:1form}.) Furthermore, $\omega_{m-1}$ involves no $dx_m$ terms and then is free of~$x_m$ by {\rm(iv)} of Lemma~\ref{LEM:form}.  Thus $\omega_{m-1}$ is a closed 1-form over $F = k(x_1,\dots,x_{m-1})$.

In Step~(5), recursively applying the above procedure over~$F$ leads to $\omega_{m-1}= dg_{m-1}+r_{m-1}$,
where $g_{m-1}\in F$ and $r_{m-1}$ is a simple form over $F$.
Set $g := a_m+g_{m-1}$ and $r := \bar{r}_m + r_{m-1}$. Thus $\omega = dg + r$, where $r$ is simple by Remark~\ref{RMK:sumsimple}. 
\end{proof}
\begin{remark}\label{RMK:rational}
	Let $\omega$ be a closed form and $\omega = dg + r$ be the decomposition computed by Algorithm~\texttt{HermiteOneForm}. By Proposition~\ref{PROP:simpleform}, $\omega$ is exact over~$K$ if and only if $r = 0$. 
It follows that the decomposition $\omega = dg + r$ is unique.  Although Step (3) of the algorithm uses the roots of the Rothstein--Trager resultant, the computation of $\bar f_i$ in Step (4) returns us to $K = k(x_1,\dots,x_m)$ with the original base field $k$.  Therefore, we have fully extended the nice features of classical Hermite reduction for rational functions. 
\end{remark}

In Section \ref{SECT:ZeilCT}, we will give a reduction-based algorithm for creative telescoping of 1-forms with one parameter.  The parameter lives in the base field, so it is important that Algorithm~\texttt{HermiteOneForm} preserves the base field.

\begin{example} \label{EXAM:Hermite}
	Let
	\[\omega = \frac{txyz-1}{x^2yz} \, dx + \frac{txyz-1}{xy^2z}\, dy + \frac{t^2xyz+xyz-1}{xyz^2}\, dz.\]
	If we let $g := \frac{1}{xyz}$ and
	\[r:= \frac{t}{x} \, dx+ \frac{t}{y} \, dy + \frac{t^2+1}{z}\, dz,\]
	then $\omega = dg+r$ is the Hermite reduction of $\omega$.
\end{example}

Example~\ref{EXAM:Telescoper} in Section~\ref{SECT:ZeilCT} will revisit this example from the point of view of creative telescoping.

\section{Integration of Closed Rational $p$-Forms}
\label{SECT:pform}

In this section, we will extend Theorem~\ref{THM:1form} to the case of rational closed $p$-forms in $m$ variables. The strategy is to reduce the problem from differential forms in $m$ variables to forms in $m-1$ or fewer variables and then proceed recursively. The result for $p$-forms differs from the case of 1-forms because the coefficients of the  logarithmic terms can be more complicated. We continue to use the notation $K = F(x_m)$ with $F = k(x_1, \ldots, x_{m-1})$, where $k$ is a field of characteristic zero.

We first give an example to show that the factors in front of the logarithms may fail to be constants when considering rational closed $p$-forms with $p > 1$.

\begin{example}
\label{EXAM:pform}
Let $K = \bC(x, y)$ and $\omega = x^{-1}y^{-1}\, dxdy$. Suppose that we can write $\omega = d(A \,dy - B \,dx)$ with
\[
A = a_0 + \sum_{i>0} \lambda_i \log a_i \quad \text{and}\quad B = b_0 + \sum_{j>0} \mu_j \log b_j,
\]
where $a_i, b_j \in \bC(x, y)$ and $\lambda_i, \mu_j$ are constants in $\bC$.  Then $x^{-1}y^{-1} = \partial_x(A) + \partial_y(B)$.  We will use residues to obtain a contradiction.  

If we view $x^{-1}y^{-1}$, $\partial_x(A)$ and $\partial_y(B)$ as rational functions in $y$ with coefficients in $\bC(x)$, then by arguments similar to the proof of Lemma~\ref{LEM:fact}, taking the residue at $y= 0$ of both sides of $x^{-1}y^{-1} = \partial_x(A) + \partial_y(B)$ leads to
\[
\frac{1}{x} =  \partial_x(\text{residue}_{y=0}(a_0)) + c \quad \text{for some constant $c\in \bC$},
\]
which is impossible by Lemma~\ref{LEM:fact}.  

However, if we no longer require that the factors in front of the logarithms be constants, then we can write
\[
\omega = d\Big(\frac{1}{y}\log (x)\, dy\Big) = d\Big({-}\frac{1}{x}\log(y)\, dx\Big).
\]
\end{example}

Now let $\omega$ be a closed $p$-form with coefficients in $K = k(x_1, \ldots, x_m) = F(x_m)$.  To explore its integration, write $\omega = \omega_{m-1} + \omega^m$ and $d = d_{m-1} + d^m$ as in Section~\ref{SECT:prelim}. Following the proof of 
Theorem 1.4.1 in~\cite{Weintraub}, we write
\[
\omega^m = dx_m(A_1 \, dx_{I_1} + \cdots + A_n \, dx_{I_n}),
\]
where $A_i \in K$ and the differentials $dx_{I_i}$ are free of $dx_m$.  By the 
Rothstein--Trager Algorithm used in the proof of Theorem~\ref{THM:1form}, we have
\[
A_i  = \partial_m(B_i) \quad \text{with}\ B_i = b_{i, 0} + \sum_{j>0} c_{i, j} \log b_{i, j},
\]
where  $b_{i, 0} \in K$, $c_{i, j} \in \overline{F}$, and $b_{i, j} \in F(c_{ i, j})[x_m]$ 
has positive degree in~$x_m$. 

Then the $(p-1)$-form $\Psi^m := B_1 \,dx_{I_1} + \cdots + B_n \,dx_{I_n}$ satisfies $\omega^m = d^m(\Psi^m)$, and we can define
\[
\widetilde{\omega} := \omega - d\hskip1pt\Psi^m = (\omega_{m-1} + \omega^m) - (d_{m-1} + d^m)(\Psi^m) = \omega_{m-1} - d_{m-1}(\Psi^m).
\]
Since $d\hskip1pt\omega =0$, we have $d\hskip1pt \widetilde{\omega}=0$. Note that $\widetilde{\omega}$ is free of $dx_m$. By {\rm(iv)} of Lemma~\ref{LEM:form}, we have $d_{m-1}(\widetilde{\omega})=0$ and $d^m(\widetilde{\omega})=0$, which implies that $\widetilde{\omega}$ is free of $x_m$, i.e., the coefficients of $\widetilde{\omega}$ are constants with respect to $x_m$. It follows that we can write $\widetilde{\omega}$ as
\[
\widetilde{\omega} = f_{\tilde I_1} \, dx_{\tilde I_1} + f_{\tilde I_2} \, dx_{\tilde I_2} + \cdots,
\]
where each $dx_{\tilde I_i}$ is free of $dx_m$ and $\partial_m(f_{\tilde I_i}) = 0$ since $d^m(\widetilde{\omega}) = 0$.

As defined, $\widetilde \omega$ appears to involve logarithms and elements of $\overline{F}$ since these appear in $\Psi^m$. We claim that the coefficients of $\widetilde{\omega}$ actually lie in $k(x_1, \ldots, x_{m-1})$.  This is proved as follows.
For $\ell \in \{1, \ldots, m-1\}$ and $i\in\{1, \ldots, n\}$, we have
\begin{align*}
\partial_\ell(B_i) &= \partial_\ell\Big(b_{i, 0} + \sum_{j>0} c_{i, j} \log b_{i, j}\Big)\\
&= \partial_\ell(b_{i, 0}) + \sum_{j>0} c_{i, j} \frac{\partial_\ell(b_{i, j})}{b_{i, j}} + \sum_{j>0} \partial_\ell(c_{i, j}) \log b_{i, j}.
\end{align*}
Thus
\begin{align*}
  \widetilde\omega  &= \omega_{m-1} - d_{m-1}(\Psi^m)  = \omega_{m-1} - d_{m-1}\bigg(\sum_{i=1}^n B_i \,dx_{I_i} \bigg) \\[-7pt]
   &  = \omega_{m-1} - \sum_{\ell=1}^{m-1}\sum_{i=1}^n \partial_\ell(B_i)\,dx_\ell dx_{I_i} \\
   & =  \omega_{m-1} - \sum_{\ell=1}^{m-1}\sum_{i=1}^n \bigg(\partial_\ell(b_{i, 0}) + \sum_{j>0} c_{i, j} \frac{\partial_\ell(b_{i, j})}{b_{i, j}} + \sum_{j>0} \partial_\ell(c_{i, j}) \log b_{i, j}\bigg)\,dx_\ell dx_{I_i} \\
   & = \widetilde{\omega}_{m-1} - \sum_{\ell=1}^{m-1}\sum_{i=1}^n \bigg(
   \sum_{j>0}c_{i,j}\frac{\partial_\ell(b_{i,j})}{b_{i,j}}+
   \sum_{j>0} \partial_\ell(c_{i, j}) \log b_{i, j}\bigg)\,dx_\ell dx_{I_i} ,
\end{align*}
where on the last line, $\widetilde{\omega}_{m-1} :=\omega_{m-1}-\sum_{\ell=1}^{m-1}\sum_{i=1}^{n}\partial_\ell(b_{i,0})\,dx_\ell dx_{I_i}$ has coefficients in $K$.
Observe that each $f_{\tilde I_i}$ is an element in $K$ minus the sum 
\[ 
\sum_{\substack {i, \ell \\ dx_\ell \wedge dx_{I_i} = \varepsilon_{i,\ell}  dx_{\tilde I_i}} }  \hskip-15pt \varepsilon_{i,\ell}  \Big(\sum_{j>0}c_{i,j}\frac{\partial_\ell(b_{i,j})}{b_{i,j}}+
\sum_{j>0} \partial_\ell(c_{i, j}) \log b_{i, j}\Big),
\] 
where $\varepsilon_{i,\ell} = \pm1$.  Since $\partial_{m}(f_{\tilde I_i}) = 0$, Lemma~\ref{LEM:fact} implies that the logarithmic part of $f_{\tilde I_i}$ vanishes.  Thus $f_{\tilde I_i}$ is an element in $K$ minus the sum 
\[
\sum_{\substack {i, \ell \\ dx_\ell \wedge dx_{I_i} = \varepsilon_{i,\ell}  dx_{\tilde I_i}} }  \hskip-15pt \varepsilon_{i,\ell}  \Big(\sum_{j>0}c_{i,j}\frac{\partial_\ell(b_{i,j})}{b_{i,j}}\Big).
\]
However, since $c_{i,j}$ and $b_{i,j}$ come from the Rothstein--Trager algorithm, the  proof of Theorem~\ref{THM:1form} shows that the sums in parentheses in the above expression are also elements of $K$.  It follows that  $f_{\tilde I_i}$ lies in $K$.
Thus $\widetilde{\omega} = \omega - d(\Psi^m)$ is free of $x_m$ and $dx_m$ and has coefficients in $k(x_1,\dots,x_{m-1})$. Furthermore, since $d = d_m$ in the $m$-variable case, we see that $\omega = \widetilde{\omega} + d_m(\Psi^m)$.

If we repeat the above process on the closed rational form $\widetilde{\omega}$ in $x_1,\dots,x_{m-1}$, we obtain $\widetilde{\widetilde{\omega}}$ and $\Psi^{m-1}$, where $\widetilde{\widetilde{\omega}}$ is a closed rational form in $x_1,\dots,x_{m-2}$ and $\Psi^{m-1}$ is similar to $\Psi^m$ except that $x_m$ and $dx_m$ no longer appear.  The result is that $\widetilde{\omega} = \widetilde{\widetilde{\omega}} + d_{m-1}(\Psi^{m-1})$, where $d_{m-1}$ is defined in Section~\ref{SECT:prelim}.  Thus
\[
\omega = \widetilde{\widetilde{\omega}} + d_{m-1}(\Psi^{m-1}) + d_m(\Psi^m).
\]
Repeating this process recursively, we must stop when only $p$ variables remain since there are no $p$-forms in fewer than $p$ variables.  Hence we have proved the following result.

\begin{thm}
\label{THM:pform}
Let $\omega$ be a closed $p$-form with coefficients in $K = k(x_1, \ldots, x_m)$. Then $\omega$ can be written in the form
\[
\omega = d_p(\Psi^p) + \cdots + d_m(\Psi^m) = d(\Psi^p + \cdots + \Psi^m),
\]
where for each $i \in \{p, \ldots, m\}$, $\Psi^i$ is free of~$dx_i,\ldots,dx_m$ and the coefficients $f_{i, \ell}$ of $\Psi^i$ are of the form
\[
f_{i, \ell} =  b_{i, \ell, 0} + \sum_{j > 0} c_{i, \ell, j} \log b_{i, \ell, j}
\]
with
\begin{align*}
b_{i, \ell, 0} &\in k(x_1, \ldots, x_{i}),\\
c_{i, \ell, j} &\in \overline{k(x_1, \ldots, x_{i-1})},\\
b_{i, \ell, j} &\in k(x_1, \ldots, x_{i-1})(c_{i, \ell, j})[x_i],\ j > 0.
\end{align*}
\end{thm}

 \begin{remark}
The above theorem says that Theorem~\ref{THM:1form} can be extended to general $p$-forms provided that the coefficients of the logarithms are allowed to be more general (specifically, algebraic over $k(x_1,\dots,x_{i-1})$ for some $i \in \{p,\dots,m\}$).
 \end{remark}

\begin{example}\label{EXAM:1/xy}
Let us apply Theorem~\ref{THM:pform} to the $2$-form $\omega = x^{-1}y^{-1}\, dxdy$ from Example~\ref{EXAM:pform}.   In the notation of the theorem with $x,y$ in place of $x_1,x_2$, we have
\[
\omega^2 = \omega = x^{-1}y^{-1}\, dxdy = dy(A_1 \, dx) \text{ with } A_1 = -x^{-1}y^{-1}.
\]
Furthermore,
\[
A_1 = \partial_y(B_1) \text{ with } B_1 = -x^{-1} \log(y),
\]
and the above procedure gives
\[
\omega = d_2(\Psi^2) = d(B_1 \, dx) = d\Big({-}\frac{1}{x}\log(y)\, dx\Big),
\]
exactly as in  Example~\ref{EXAM:pform}. Switching $x$ and $y$ gives the other representation of $\omega$ in the example.
\end{example}

\begin{example}
\label{EXAM:2formalgebraic}
Algebraic functions can also appear in front of the logarithms.  Consider the closed 2-form
\[
\omega = \frac{1}{x^5 + x - y^2}\,dxdy
\]
Using $x,y$ in place of $x_1,x_2$, Theorem~\ref{THM:pform} leads to
\begin{align*}
\omega = d\Psi^2 &= d\Big(\!\Big(\frac{1}{2\sqrt{x^5+x}}\log(y-\sqrt{x^5+x}) - \frac{1}{2\sqrt{x^5+x}}\log(y+\sqrt{x^5+x})\Big)\hskip1pt dx\Big)\\
&= d\Big(\frac{1}{2\sqrt{x^5+x}}\log \frac{y-\sqrt{x^5+x}}{y+\sqrt{x^5+x}}\,dx\Big).
\end{align*}
Switching $x$ and $y$ gives an answer that uses the five roots $x_1(y),\dots,x_5(y)$ of $x^5 + x - y^2 = 0$ regarded as algebraic functions of $y$.  The procedure of the theorem gives
\[
\omega = d\Big(\!\Big(\frac{1}{5x_1(y)^4+ 1}\log(x - x_1(y)) + \cdots + \frac{1}{5x_5(y)^4+ 1}\log(x - x_5(y)) \Big)\hskip1pt dy\Big).
\]
These formulas make it clear that arbitrarily complicated algebraic functions can appear as the coefficients of the logarithms when dealing with $p$-forms for $p \ge 2$.
\end{example}

\begin{example}
\label{EXAM:2form}
Consider the closed 2-form
\[
\omega = \big( \frac{1}{z^2-x} + \frac{1}{x y}\big)\, dxdy + \frac{1}{z^2-x} \, dydz + \frac{y-2yz}{(z^2-x)^2}\, dxdz.
\]
In the notation of Theorem~\ref{THM:pform} with $x,y,z$ in place of $x_1,x_2,x_3$ respectively, we have
$\omega = \big( 1/(z^2-x) + 1/(xy) \big)\, dxdy + dz(A_1 \, dx + A_2 \, dy)$  with
\[
A_1 := -\frac{y-2yz}{(z^2-x)^2}  \text{ and } A_2 := -\frac{1}{z^2-x}.
\]
Furthermore, $A_1 = \partial_z(B_1)$ with
\[
B_1 := \frac{\sqrt{x}y}{4x^2} \log\frac{z-\sqrt{x}}{z+\sqrt{x}} + \frac{ yz - 2xy}{2x(z^2-x)},
\]
and~$A_2 = \partial_z(B_2)$ with
\[
B_2 := -\frac{\sqrt{x}}{2x} \log\frac{z-\sqrt{x}}{z+\sqrt{x}} .
\]
Hence, $\Psi^3 := B_1 \, dx + B_2 \, dy$ and
\[
\widetilde \omega := \omega - d(\Psi^3) = \bigg( \frac{1}{z^2-x} + \frac{1}{x y} + \partial_y(B_1)  - \partial_x(B_2) \bigg) \, dxdy = \frac{1}{xy}\, dxdy.
\]
Example~\ref{EXAM:1/xy} shows that $\widetilde{\omega} = d(\Psi^2)$ with $\Psi^2 := {-}\frac1x \log(y) \,dx$, so  $\omega =  \widetilde{\omega} + d(\Psi^3)$ is given by
\begin{align*}
\omega &= d(\Psi^2) + d(\Psi^3) =d(\Psi^2 +\Psi^3)\\
&= d\biggl(\! \Bigl({-} \frac{1}{x}\log(y) + \frac{\sqrt{x}y}{4x^2} \log\frac{z-\sqrt{x}}{z+\sqrt{x}} + \frac{ yz - 2xy}{2x(z^2-x)} \Bigr)\, dx -\frac{\sqrt{x}}{2x} \log\frac{z-\sqrt{x}}{z+\sqrt{x}} \, dy  \biggr).
\end{align*}
Notice that the coefficients of the logs in $\Psi^3$ involve $x$ and $y$, while the coefficient of the log in $\Psi^2$ involves only $x$, exactly as described in Theorem~\ref{THM:pform}.  

However, if we take $z,y,x$ in place of~$x_1,x_2,x_3$ respectively, then the algorithm of Theorem~\ref{THM:pform} leads to the following simpler expression:
\[
    \omega = d\biggl(\!\Bigl(-\log(z^2-x)+\frac{1}{y}\log(x)\Bigr)\, dy + \frac{y-2yz}{z^2-x} \, dz\biggr).
\]
The 1-form on the right involves fewer logarithms than the above 1-form $\Psi^2 + \Psi^3$.
\end{example}

\begin{example}
\label{EXAM:2form2}
Removing $1/(xy)\,dxdy$ from $\omega$ in Example~\ref{EXAM:2form} gives the closed 2-form
\[
\omega' =  \frac{1}{z^2-x}\, dxdy + \frac{1}{z^2-x} \, dydz + \frac{y-2yz}{(z^2-x)^2}\, dxdz,
\]
and the computation of Example \ref{EXAM:2form} shows that
\[
\omega' = d\Psi^3 = d\biggl(\! \Bigl(\frac{\sqrt{x}y}{4x^2} \log\frac{z-\sqrt{x}}{z+\sqrt{x}} + \frac{ yz - 2xy}{2x(z^2-x)} \Bigr)\, dx -\frac{\sqrt{x}}{2x} \log\frac{z-\sqrt{x}}{z+\sqrt{x}} \, dy  \biggr).
\]
Here, $\Psi^3$ involves the rational function $\dfrac{ yz - 2xy}{2x(z^2-x)}$, whose denominator contains the factor $x$ not present in the denominators of $\omega'$.  
\end{example}

Example \ref{EXAM:2form2} shows that the algorithm of Theorem~\ref{THM:pform} can introduce new denominators, and Examples~\ref{EXAM:2formalgebraic} and  \ref{EXAM:2form} show that the output of the algorithm depends on the ordering of the variables and that different orderings can lead to different numbers of logarithms.  These examples raise two natural questions about what happens when we write a closed rational $p$-form $\omega$ as $\omega = d\Psi$:
\begin{enumerate}
    \item Can $\Psi$ be chosen so as to not introduce any new denominators?
    \item Can $\Psi$ be chosen so as to minimize the number of logarithms involved?
\end{enumerate}
Example \ref{EXAM:ratint} below will illustrate the challenge posed by question (2).

\section{Picard's Problem and Griffiths--Dwork Reduction}
\label{SECT:Picard}

Let $k$ be a field of characteristic zero and $k(x, y, z)$ be the field of rational functions in $x, y$ and $z$ over $k$. Let $\partial_x, \partial_y, \partial_z$ denote the usual derivations on $k(x, y, z)$ with respect to $x, y, z$, respectively. In his books~\cite[Vol.\ II, pp.~475--479]{picard1897} and~\cite[pp.~209--215]{picard1931}, Picard posed the following problem.

\begin{PP}
 Given $f\in k(x, y, z)$, decide whether there exist $u, v, w \in k(x, y, z)$ such that
\begin{equation}\label{EQ:Picard}
f = \partial_x(u) + \partial_y(v) + \partial_z(w).
\end{equation}
\end{PP}

In this section, we will study an explicit example and reformulate Picard's problem in terms of differential forms. In the case of $m$ variables, we will then discuss how recent work on  the 
Griffiths--Dwork order of pole reduction relates to Picard's problem. Note that the discrete analogue of Picard's problem has been recently solved in~\cite{ChenDuFang2025}.

We begin with some definitions.
Set $K=k(x, y)$ and let $\overline{K}$ be the algebraic closure of $K$. Any rational function $f=P/Q \in K(z)$ can be decomposed into
\[
  f = p + \sum_{i=1}^m \sum_{j=1}^{n_i}\frac{\alpha_{i,j}}{(z-\beta_i)^j},
\]
where $p\in {K}[z]$, the $\beta_i$'s are the roots in $\overline{K}$ of the denominator $Q$, and $\alpha_{i,j}\in \overline{K}$. We call the element $\alpha_{i,1}$ the \textit{$z$-residue} of $f$ at $\beta_i$.

A rational function $f\in k(x, y, z)$ is said to be \textit{rationally integrable} with respect to $x, y, z$ if $f=\partial_x(u)+\partial_y(v)+\partial_z(w)$ for some $u, v, w\in k(x, y, z)$. An algebraic function $\alpha\in
\overline{k(x, y)}$ is said to be \textit{algebraically integrable} with respect to $x, y$ if $\alpha =\partial_x(\beta) + \partial_y(\gamma)$ for some $\beta, \gamma\in \overline{k(x, y)}$. We recall Lemma 4 from~\cite{ChenKauersSinger2012} on the relation between those two notions of being integrable.

\begin{lem}
\label{LEM:ratsol}
Let $f\in k(x, y, z)$. Then $f$ is rationally integrable with respect to $x, y, z$ if and only if all $z$-residues of $f$ are algebraically integrable with respect to~$x, y$.
\end{lem}

We now show that some rational functions are not rationally integrable.

\begin{example}
\label{EXAM:notint}
For a fixed $n\in \NX$, consider the rational function
\[
f = \frac{1}{x^n + y^n + z^n}.
\]
We will prove that $f$ is rationally integrable with respect to $x, y, z$ if and only if $n\neq 3$. It is straightforward to check that when $n\neq 3$ we have
\[
f = \partial_x\left(\frac{(3-n)^{-1}x}{x^n + y^n + z^n}\right) +  \partial_y\left(\frac{(3-n)^{-1}y}{x^n + y^n + z^n}\right) + \partial_z\left(\frac{(3-n)^{-1}z}{x^n + y^n + z^n}\right).
\]
So it remains to show that $f$ is not rationally integrable with respect to $x, y, z$ when $n=3$. A root (in $\overline{k(x, y)}$) of the denominator of $f = 1/(x^3+y^3+z^3)$ is $\lambda (x^3 + y^3)^{1/3}$ with $\lambda^3=-1$. Then the residue of $f$ at this root is
\[
\frac{1}{3\lambda^{2}(x^3 + y^3)^{2/3}}.
\]
By Lemma~\ref{LEM:ratsol}, it suffices to show that the algebraic function
\[
\alpha(x, y)= \frac{1}{(x^3 + y^3)^{2/3}}
\]
is not algebraically integrable. Consider the birational transformation
\[\bar{x} = x \quad \text{and} \quad \bar{y} = y/x.\]
Then one computes that
\begin{equation}\label{EQ:alphabaralpha}
\alpha(x, y) \,dx  dy = \bar{\alpha}(\bar x, \bar y)\,d\bar x  d \bar y, \quad \text{where}\ \bar \alpha (\bar x, \bar y)= \frac{1}{\bar{x}(1+\bar{y}^3)^{2/3}}.
\end{equation}
We will first show that $\bar \alpha$ is not algebraically integrable with respect to $\bar{x}, \bar{y}$.  

Seeking a contradiction, assume that $\bar \alpha$ is algebraically integrable, so that $\bar \alpha = \partial_{\bar x}(\bar \beta) + \partial_{\bar y}(\bar \gamma)$ for some $\bar \beta, \bar \gamma\in \overline{k(\bar x, \bar y)}$.  We cannot use Lemma~\ref{LEM:ratsol} directly, so we instead adapt the proof of~\cite[Lemma 4]{ChenKauersSinger2012}. Note that
\[
 \bar \beta, \bar \gamma\in\overline{k(\bar x, \bar y)} = \overline{\overline{k(\bar y)} (\bar x)} \subseteq \bigcup_{n=1}^\infty \overline{k(\bar y)}(({\bar x}^{1/n})),
\]
where on the right we have the field of Puiseux series in $\bar x$ with coefficients in $\overline{k(\bar y)}$.  For $\bar \delta \in \overline{k(\bar x, \bar y)}$, define $\text{residue}_{\bar x}(\bar \delta,0)$ to be the coefficient of ${\bar x}^{-1}$ in the Puiseux series of $\bar\delta$ in fractional powers of $\bar x$.  Note that $\text{residue}_{\bar x}(\bar\delta,0) \in \overline{k(\bar y)}$, and similar to~\cite[Prop.~3]{ChenKauersSinger2012}, one can easily show that
\[
\text{residue}_{\bar x}(\partial_{\bar x}(\bar\delta),0) = 0 \quad \text{and} \quad \text{residue}_{\bar x}(\partial_{\bar y}(\bar\delta),0) = \partial_{\bar y}(\text{residue}_{\bar x}(\bar\delta,0)).
\]
From $\bar \alpha = 1/(\bar x(1+\bar{y}^3)^{2/3}) = \partial_{\bar x}(\bar \beta) + \partial_{\bar y}(\bar \gamma)$, we obtain
\begin{align*}
\frac{1}{(1+\bar{y}^3)^{2/3}} &= \text{residue}_{\bar x}(\bar \alpha,0) = \text{residue}_{\bar x}(\partial_{\bar x}(\bar \beta) + \partial_{\bar y}(\bar \gamma),0)\\ &= \text{residue}_{\bar x}(\partial_{\bar y}(\bar \gamma),0) = \partial_{\bar y}(\text{residue}_{\bar x}(\bar\gamma,0)),
\end{align*}
which proves that
\[
\beta(\bar{y}) = \frac{1}{(1+\bar{y}^3)^{2/3}} = \partial_{\bar y}(\text{residue}_{\bar x}(\bar\gamma,0))
\]
is algebraically integrable with respect to $\bar y$ since $\text{residue}_{\bar x}(\bar\gamma,0) \in \overline{k(\bar y)}$.

We next show that the simple radical $\beta(\bar{y})$ is not algebraically integrable with respect to $\bar{y}$. In fact, if $\beta(\bar{y})$ is algebraically integrable, then
\[ \beta(\bar{y})=D_{\bar{y}} (\gamma(\bar y)), \]
where $\gamma(\bar y)$ is algebraic over $k(\bar y)$ and $D_{\bar y} = d/d{\bar y}$.   By~\cite[Lemma 5.1.1]{Bronstein}, we can assume that $\gamma(\bar y) \in k(\bar y, \beta(\bar y))$. 
Using $D_{\bar y}(\beta) = \frac{-2\bar{y}^2}{1+\bar{y}^3} \beta$ and the basis $\{1,\beta(\bar y),\beta(\bar y)^2\}$ of $k(\bar y, \beta(\bar y))$ over $k(\bar y)$, one can show that $\gamma(\bar y) = r\beta(\bar y) + c$, where  $r \in k(\bar{y})$ and $c \in k$.  Thus
\[ 
\beta(\bar{y})=D_{\bar{y}} \bigg(\frac{r}{(1+\bar{y}^3)^{2/3}}\bigg), 
\]
where $r \in k(\bar{y})$.  It follows that
\begin{equation}\label{EQ:nequals3}
D_{\bar{y}}(r)-\frac{2\bar{y}^2}{1+\bar{y}^3}\cdot r =1.
\end{equation}

A straightforward application of the Risch Algorithm from~\cite[Ch.\ 6]{Bronstein} to \eqref{EQ:nequals3} shows that the rational solution $r$ is in fact a polynomial of degree at most 2 (see Def.\ 6.1.1, Thm.\ 6.1.2, and Cor.\ 6.3.1 of~\cite{Bronstein}).  An easy computation now leads to a contradiction, so that \eqref{EQ:nequals3} has no rational solution. As we have seen, this implies that $\bar \alpha$ is not algebraically integrable.

We claim that the same is true for $\alpha$.   Here is a classical argument using Picard and Simart~\cite{picard1897} (see Remark~\ref{RMK:algebraicallyexact} below for a modern proof). Suppose that $\alpha = \partial_x(\beta) + \partial_y(\gamma)$ for $\beta,\gamma \in \overline{k(x, y)}$.  In terms of integrals, this means that
\[
\iint \alpha\,dxdy = \iint  \partial_x(\beta) + \partial_y(\gamma)\,dxdy.
\]
By~\cite[Vol.\ II, pp.~160--161]{picard1897}, integrals such as the one on the right are preserved under birational transformations when $\beta$ and $\gamma$ are rational functions.  The proof adapts easily to the case of algebraic functions, so that for the birational transformation given above, we have
\[
 \iint  \partial_x(\beta) + \partial_y(\gamma)\,dxdy =  \iint  \partial_{\bar x}(\bar\beta) + \partial_{\bar y}(\bar\gamma)\,d\bar x d\bar y
\]
for algebraic functions $\bar\beta$ and $\bar\gamma$.  Combining this with \eqref{EQ:alphabaralpha}, we conclude that $\bar \alpha$ is algebraically integrable, which contradicts what we have proved above.

It follows that $\alpha$ is not algebraically integrable.  As noted above, Lemma~\ref{LEM:ratsol} then implies that $f = 1/(x^3+y^3+z^3)$ is not rationally integrable.
\end{example}

The $2$-form $(\partial_x(\beta) + \partial_y(\gamma))\,dxdy$ appearing in Example~\ref{EXAM:notint} suggests that the $3$-form $(\partial_x(u) + \partial_y(v) +\partial_z(w))\,dxdydz$ may be relevant to Picard's problem.  In fact,
\[
d(u \,dy  dz - v \,dx  dz + w \,dx dy) = (\partial_x(u) + \partial_y(v) +\partial_z(w))\,dxdydz,
\]
which implies that Picard's problem \eqref{EQ:Picard} is equivalent to saying that the $3$-form $f(x, y, z)\, dx dy  dz$ is exact over~$k(x,y,z)$, i.e.,
\[
f \,dx dy  dz = d(u \,dy  dz - v \,dx  dz + w \,dx dy),\quad\text{where}\ u,v,w \in k(x,y,z).
\]

\begin{remark}\label{RMK:algebraicallyexact}
Similar formulas also apply to algebraic functions in two variables, so that $\alpha(x,y)$ from
Example~\ref{EXAM:notint} is algebraically integrable if and only if the $2$-form $\alpha(x,y)\,dxdy$ is the exterior derivative of an algebraic $1$-form.  Since the exterior derivative is compatible with pullbacks by Remark~\ref{RMK:pullback}, we get an immediate proof that in \eqref{EQ:alphabaralpha}, $\alpha(x,y)$ is not  algebraically integrable if and only if the same is true for $\bar\alpha(\bar x,\bar y)$.
\end{remark}

We can also relate Example~\ref{EXAM:notint} to Theorem~\ref{THM:pform} as follows.

\begin{example}
\label{EXAM:3Form}
By the previous example, $f = 1/(x^3 + y^3 + z^3)$ is not  rationally integrable with respect to $x, y, z$, so that $f\, dxdydz$ is not the exterior derivative of a rational $2$-form.  However, Theorem~\ref{THM:pform} shows how to write $f\, dxdydz$ as  the exterior derivative of a rational and logarithmic $2$-form as follows.

Let $\alpha = -(x^3 + y^3)^{1/3}$. Then the three roots in $\overline{\QX(x, y)}$ of $z^3 + x^3 + y^3 = 0$ are $z = \alpha, \alpha \omega, \alpha \omega^2$, where $\omega$ is a primitive cube root of unity. Then
\[
f = \frac{1/(3\alpha^2)}{z-\alpha} +  \frac{\omega/(3\alpha^2)}{z-\alpha\omega} + \frac{\omega^2/(3\alpha^2)}{z-\alpha\omega^2},
\]
which implies that $f = \partial_z(g)$ with
\[
g = \frac{1}{3\alpha^2}\log(z-\alpha) + \frac{\omega}{3\alpha^2}\log(z-\alpha\omega) + \frac{\omega^2}{3\alpha^2}\log(z-\alpha\omega^2) .
\]
If follows easily that $f \,dx  dy dz = d(g\, dx  dy)$.
\end{example}

\begin{conj}
\label{CONJ:WhenRatInt}
For a fixed $n \in \mathbb{N}$, the rational function
\[f = \frac{1}{x_1^n+\cdots+x_m^n} \]
is rationally integrable with respect to $x_1,\ldots,x_m$ with $m \geq 4$ if and only if $n \neq m$.
\end{conj}

The general case of Picard's problem asks when $f \in K = k(x_1, \dots, x_m)$ can be written in the form
\begin{equation}
\label{EQ:genratint}
f = \partial_{1}(u_1) + \cdots + \partial_{m}(u_m),\quad\text{where}\ u_1,\dots,u_m \in K.
\end{equation}
Since
\[
d\Big(\sum_{i=0}^m (-1)^{i-1} u_i \, dx_1 \cdots \widehat{dx_i} \cdots dx_m\Big) = \Big(\sum_{i=1}^m \partial_{i}(u_i)\Big)\,dx_1 \cdots  dx_m
\]
(the $\widehat{\hphantom{dx}}$ over~$dx_i$ indicates that it is omitted),
Picard's problem is thus equivalent to deciding when the $m$-form
\[
\omega = f\, dx_1 \cdots  dx_m, \quad\text{where}\  f \in K,
\]
is rationally integrable, i.e., is the exterior derivative of a rational $(m-1)$-form.

We next discuss a method relevant to Picard's problem.  This method, sometimes called \textit{Griffiths--Dwork order of pole reduction}, has a long history that goes back to Hermite~\cite{Hermite1872} in 1872, Dwork~\cite{Dwork1962,Dwork1964} in 1962 and 1964, Griffiths~\cite{Griffiths1969} in 1969, and Dimca~\cite{Dimca1991} in 1991.  Our exposition will focus on the more recent work by Bostan-Lairez-Salvy~\cite{BostanPierreBruno2013} in 2013 and especially Lairez~\cite{Lairez2016} in 2016.  

Our first task is to translate Picard's problem into the homogeneous setting used in 
\cite{BostanPierreBruno2013,Lairez2016}.  Homogeneous elements $R \in k(x_0,x_1,\dots,x_m)$ of degree $-m-1$ play a special role because the $m$-form $R\, \Omega$, where
\[
\Omega = \sum_{i=0}^m (-1)^i x_i\,dx_0 \cdots \widehat{dx_i} \cdots dx_m,
\]
gives a rational $m$-form on $\mathbb{P}^m$, and all rational $m$-forms on $\mathbb{P}^m$ arise in this way.  This observation (see~\cite[Cor.\ 2.11]{Griffiths1969}) is the starting point of the Griffiths--Dwork method in ~\cite{Griffiths1969}.  We will instead follow the approach of \cite{Lairez2016}, which for $R$ as above works with the 
$(m+1)$-form $R\,dx_0dx_1\cdots dx_m$ of degree $0$.

We will need the following notation.  First, let $Q$ be a squarefree homogeneous polynomial in $k[x_0,x_1,\dots,x_m]$, and let 
\[
A_Q = k[x_0,x_1,\dots,x_m]_Q = k[x_0,x_1,\dots,x_m,1/Q]
\]
be the localization of $k[x_0,x_1,\dots,x_m]$ at $Q$.  Next, for $a \in \mathbb{Z}$, 
define the homogenization operator $h_a$ on $k(x_1,\dots,x_m)$ as follows: for $g \in k(x_1,\dots,x_m)$, set
\[
h_a(g) = x_0^{a} \, g(x_1/x_0,\dots,x_m/x_0).
\]
Note that $h_a(g)  \in k(x_0,x_1,\dots,x_m)$ is homogeneous of degree $a$.  

\begin{lem}
\label{LEM:HomogeneousPicard}
Let $f \in k(x_1,\dots,x_m)$ and let $G := h_{-m-1}(f) \in k(x_0,x_1,\dots,x_m)$ be its homogenization of degree $-m-1$.  Then 
$f$ is rationally integrable if and only if there exists $Q$ homogeneous and squarefree such that $G \in \sum_{i=0}^m \partial_i(A_Q)$.
\end{lem}

\begin{proof}
First suppose that $f$ is rationally integrable, so that $f = \sum_{i=1}^m \partial_i(u_i)$ where $u_i \in k(x_1,\dots,x_m)$.  Let $P$ be the product of the distinct irreducible factors that appear in the denominators of the $u_i$, and let $Q = x_0 h_{\deg{P}}(P)$. Note that $Q \in k[x_0,\dots,x_m]$ is homogeneous and squarefree.

Let $\tilde u_i = h_{-m}(u_i)$ and note that $\tilde u_i \in A_Q$ since every irreducible factor appearing in the denominator of $\tilde u_i$ divides $Q$.  Furthermore, observe that for $g \in k(x_1,\dots,x_m)$ and $i \in \{1,\dots,m\}$, we have
\begin{align*}
\partial_i(h_{-m}(g)) &= \partial_i( x_0^{-m} g(x_1/x_0,\dots,x_m/x_0))=  x_0^{-m}  (\partial_i(g))(x_1/x_0,\dots,x_m/x_0)\cdot 1/x_0\\ &= x_0^{-m-1}  (\partial_i(g))(x_1/x_0,\dots,x_m/x_0) = h_{-m-1}(\partial_i(g)), 
\end{align*}
where the second equality on the first line uses the chain rule. Thus
\[
G = h_{-m-1}(f) = \sum_{i=1}^m h_{-m-1}(\partial_i(u_i)) = \sum_{i=1}^m \partial_i(h_{-m}(u_i)) =  \sum_{i=1}^m \partial_i(\tilde u_i).
\]
It follows that $G\in \sum_{i=0}^m \partial_i(A_Q)$ since $\tilde u_i \in A_Q$ as noted above.

Conversely, suppose that $G = \sum_{i=0}^m \partial_i(\tilde u_i)$ for some $\tilde u_0, \tilde u_1,\dots,\tilde u_m \in A_Q$.  Every element of $A_Q$ is a sum of homogeneous components, and since $\partial_i$ lowers the degree by 1,  we can assume that each $\tilde u_i$ is homogeneous of degree $-m$.

We dehomogenize by setting $x_0 = 1$.  Note that this operation commutes with $\partial_i$ for $i \in \{1,\dots,m\}$ (but not for $i =0$).  Setting $u_i = (\tilde u_i)_{|x_0=1} \in k(x_1,\dots,x_m)$, we obtain
\begin{equation}
\label{EQ:dehomog}
\begin{aligned}
f = h_{-m-1}(f)_{|x_0=1} = G_{|x_0=1} &= (\partial_0(\tilde u_0))_{|x_0=1} + \sum_{i=1}^m (\partial_i(\tilde u_i))_{|x_0=1}\\
&= (\partial_0(\tilde u_0))_{|x_0=1} + \sum_{i=1}^m \partial_i(\,\underbrace{(\tilde u_i)_{|x_0=1}}_{\text{\normalsize $u_i$}}\,)\\ 
\end{aligned}
\end{equation}
To untangle $(\partial_0(\tilde u_0))_{|x_0=1}$, remember that $\tilde u_0$ is homogeneous of degree $-m$, so that by Euler's formula, 
\[
-m \tilde u_0 = \sum_{i=0}^m x_i \partial_i( \tilde u_0) = x_0\partial_0( \tilde u_0) + \sum_{i=1}^m x_i \partial_i( \tilde u_0).
\]
We can rewrite this as
\[
x_0 \partial_0( \tilde u_0) = -m\tilde u_0 - \sum_{i=1}^m x_i \partial_i( \tilde u_0)
=-\sum_{i=1}^m (\tilde u_0 + x_i \partial_i( \tilde u_0)) = -\sum_{i=1}^m \partial_i(x_i  \tilde u_0).
\]
Setting $x_0 = 1$ and using $u_0 = (\tilde u_0)_{|x_0=1}$, we get
\[
(\partial_0(\tilde u_0))_{|x_0=1} = -\sum_{i=1}^m \partial_i(x_i  u_0),
\]
and then combining this with~\eqref{EQ:dehomog} gives
\[
f =  \sum_{i=1}^m \partial_i(v_i), \quad v_i = u_i - x_i u_0 \in k(x_1,\dots,x_m).
\]
Thus $f$ is rationally integrable.
\end{proof}

In the setting of the lemma, note that $G \in \sum_{i=0}^m \partial_i(A_Q)$ implies that $G \in A_Q$.  In general, given any $R \in A_Q$, there is an algorithm for determining whether or not $R \in \sum_{i=0}^m \partial_i(A_Q)$.  This follows from the results of~\cite{Lairez2016}.  Specifically, Theorem 4 of~\cite{Lairez2016} shows that $G \in \sum_{i=0}^m \partial_i(A_Q)$ is equivalent to a condition that can be checked by standard methods of symbolic computation, provided that a certain integer $r$ is sufficiently large.  Then Section 5 of~\cite{Lairez2016} shows how to compute a suitable $r$.


On the surface, combining this with Lemma~\ref{LEM:HomogeneousPicard} seems to give an algorithm for solving Picard's problem.  
The difficulty lies in the phrase ``there exists'' that appears in the equivalence proved in the Lemma~\ref{LEM:HomogeneousPicard}.  Namely, given $f$, we can write $G = h_{-m-1}(f) = a/Q^\ell$ for some squarefree homogeneous polynomial $Q$, and we can test whether $G \in \sum_{i=0}^m \partial_i(A_Q)$ using~\cite{Lairez2016} as above. If the answer is ``yes'', then $f$ is rationally integrable by the lemma.  But if the answer is ``no'', then we can't guarantee that $f$ is not rationally integrable, since in a sum $f = \sum_{i=1}^m \partial_i(u_i)$, the denominators of the $u_i$ might involve irreducible factors not accounted for when we wrote $G = a/Q^\ell$.  In fact, the $Q$ constructed in the proof of Lemma~\ref{LEM:HomogeneousPicard} makes explicit uses of the denominators of the $u_i$.

On way to think about this is that if we write $G = h_{-m-1}(f) = a/Q^\ell$, then determining when $G \in \sum_{i=0}^m \partial_i(A_Q)$ can be thought of as the \textit{constrained version of Picard's problem} for $f$, where the solution does not introduce any new poles.  The paper~\cite{Lairez2016} provides an algorithm to decide this, and, as pointed out by one of the reviewers, one can also decide the constrained Picard's problem using $D$-module algorithms for computing the de Rham cohomology of a hypersurface complement (see, for example,~\cite{OakuTakayama1999}). 
 


Another approach to Picard's problem (the original unconstrained version) comes from Theorem~\ref{THM:pform}.  As noted above, Picard's problem for $f \in k(x_1,\dots,x_m)$ is equivalent to asking if the $m$ form $\omega = f\,dx_1\cdots dx_m$ is the exterior derivative of a rational $(m-1)$-form.  By Theorem~\ref{THM:pform}, we can write $\omega = d\Psi^m$, where $\Psi^m$ is a $(m-1)$-form that may involve logarithms.  If it doesn't, then $f$ is rationally integrable.  But if it does, it is possible that $f$ could still be rationally integrable, as shown by the following simple example. 

\begin{example}\label{EXAM:ratint}
Consider $\omega = 1/(x+y)\,dxdy$.  Applying Theorem~\ref{THM:pform} gives
\[
\omega = d({-}\log(x+y)\,dx),
\]
and reversing the order of the variables gives
\[
\omega = d(\log(x+y)\,dy).
\]
Logarithms appear to be inescapable.  Yet $1/(x+y)$ is rationally integrable because
\[
 \omega =  \frac{1}{x+y} \, dxdy = d \Big(\frac{x}{x+y} \, dx - \frac{y}{x+y} \, dy\Big).
 \]
\end{example}

In general, when we write $f \,dx_1\cdots dx_m = d\Psi^m$, there are many choices for $\Psi^m$, as illustrated here and in the examples given in Section~\ref{SECT:pform}.  At the end of Section~\ref{SECT:pform}, we asked if it is possible to construct a choice for $\Psi^m$ that involves as few logarithms as possible.  A positive answer to this question would provide a method for determining when $f$ is rationally integrable.



The upshot is that the problem of determining the rational integrability of rational functions in $k(x_1,\dots,x_m)$ (the unconstrained Picard's problem) remains open.

\section{Liouville's Theorem for Differential 1-forms }
\label{SECT:Liouville}

We will extend Liouville's classical theorem to
integrals of differential $1$-forms in~$\mathcal{M} = \mathcal{U}\,dx_1\oplus\cdots\oplus \mathcal{U}\,dx_m$. Recall that $K = k(x_1,\ldots,x_m)$ and $K \subseteq  \mathcal{U}$ is the universal differential extension from Section~\ref{SECT:prelim}. Throughout this section, all fields are assumed to be of characteristic zero.

\begin{defin} 
Let $E\subseteq \mathcal{U}$ be a differential extension of $K$ and $t\in E$.
\begin{enumerate}
\item $t$ is a \textit{logarithm} over $K$ if there exists nonzero $b\in K$ such that $\partial_i t=\partial_i b/b$ for each $i \in\{1,\ldots,m\}$.
\item $t$ is an \textit{exponential} over $K$ if there exists $b\in K$ such that $\partial_i t/t = \partial_i b$ for each $i \in\{1,\ldots,m\}$.
\item $E$ is an \textit{elementary extension} of $K$ if there exist $t_1, \ldots, t_n\in E$ such that $E=K(t_1,\ldots,t_n)$ and for each $j \in \{1,\ldots,n\}$, $t_j$ is either algebraic or a logarithm or an exponential over $K(t_1,\ldots,t_{j-1})$.
\end{enumerate}
\end{defin}

The following theorem is Risch's strong version of the classical Liouville's theorem~(see~\cite{Rosenlicht1968,Risch1969} and Section~5.5 in~\cite{Bronstein} for its proof).

\begin{thm}[Liouville's theorem]\label{THM:liouville}
	Let $(F,D)$ be a differential field and $f\in F$. If there exists an elementary extension $E$ of $F$ and $g\in E$ such that $D(g)=f$, then there are $v\in K$, $c_1,\ldots,c_n \in \overline{\Const_D(F)}$ and nonzero $u_1,\ldots,u_n \in F(c_1,\ldots,c_n)$ such that
	\[f = D(v) + \sum_{i=1}^{n} c_i \frac{D(u_i)}{u_i}.\]
\end{thm}

The following theorem for 1-forms was first proved by Caviness and Rothstein in~\cite{Caviness1975} and later proved by Rosenlicht~\cite{Rosenlicht1976} based on Ax's results in~\cite{Ax1971}. Both proofs proceed by induction on the number of generators of the field extension. Instead, we now present a new proof by induction on the number of derivations, based on the classical Liouville's theorem and the recursive technique used in the proof of Theorem~\ref{THM:1form}. 

\begin{thm}[Liouville's theorem for 1-forms]\label{THM:Liouville1form}
	Let $F \subseteq \mathcal{U}$ be a differential extension of~$K$  and $C := \Const_{\partial_1,\dots,\partial_m}(F)$. Let $\omega=f_1 \, dx_1+\cdots+f_m \, dx_m \in \mathcal{M}$ be a closed 1-form with $f_i \in F$. If there exist an elementary
extension $E$ of $F$ and $g\in E$ such that $\omega =dg$, then there are $a\in F$, $c_1,\ldots,c_n \in \overline{C}$, and nonzero $b_1,\ldots,b_n\in F(c_1,\ldots,c_n)$ such that
	\[g = a+\sum_{j}c_j \log b_j.\]
\end{thm}

\begin{proof}
	Denote $C_{m}:=\Const_{\partial_m}(F)$. Since $f_m$ has an elementary integral over $F$ with respect to $\partial_m $, Theorem~\ref{THM:liouville} implies that there exist $a_m \in F,\ c_{m,1},\ldots,c_{m,n} \in \overline{C_m}$ and nonzero $b_{m,1},\ldots,b_{m,n}\in F(c_{m,1},\ldots,c_{m,n})$ such that
	\[f_m=\partial_m(g_{m}) \quad \text{with} \quad g_m=a_{m}+\sum_{j=1}^{n}c_{m,j}\log{b_{m,j}} \]
	where $\partial_m(b_{m,j})\neq 0$ for each $j$.
By Corollary 3.3.1 from~\cite{Bronstein}, we have $\Const_{\partial_m}(\overline{F}) = \overline{C_m}$, so $\overline{F}' := \{\partial_m(f)\ | \ f\in \overline{F}\}$ is a linear subspace of $\overline{F}$ over $\overline{C_m}$. For every $f\in \overline{F}$, let $\widehat{f}$ be the image of $f$ in the quotient space $\overline{F}/\overline{F}'$. We can always select maximally many of the $\widehat{\frac{\partial_m(b_{m,j})}{b_{m,j}}}$'s so that they are linearly independent over $\overline{C_m}$ and then rearrange the logarithmic part of $g_m$. After this process, we can always assume that the $\widehat{\frac{\partial_m(b_{m,j})}{b_{m,j}}}$'s appearing in $g_m$ are linearly independent over $\overline{C_m}$, though $a_m$ may now be in $\overline{F}$. Since $\omega$ is closed, we have $\partial_m(f_i) = \partial_i(f_m)$ for all $i \in\{1,\ldots, m\}$.  Then \eqref{EQ:fm} and \eqref{EQ:partialfm} imply that
\[
\partial_m(f_i) =\partial_m(\partial_i(a_m))+\sum_{j=1}^{n}\partial_i(c_{m,j})\frac{\partial_m(b_{m,j})}{b_{m,j}}+\partial_m\Big(\sum_{j=1}^{n}c_{m,j}\frac{\partial_i(b_{m,j})}{b_{m,j}}\Big).
\]
Taking the image of both sides in the quotient space $\overline{F}/\overline{F}'$ gives
	\[\sum_{j=1}^{n}\partial_i(c_{m,j})\widehat{\frac{\partial_m(b_{m,j})}{b_{m,j}}}=0,\]
where $\partial_i(c_{m,j}) \in \overline{C_m}$ because $c_{m,j}\in \overline{C_m}$. Since the $\widehat{\frac{\partial_m(b_{m,j})}{b_{m,j}}}$'s are linearly independent over $\overline{C_m}$, $\partial_i(c_{m,j})=0$ for each $j$, which implies that $c_{m,j} \in \overline{C_i}$. Since this holds for all $i \in\{1,\ldots,m-1\}$, we see that $c_{m,j}\in \overline{C} $.
	
	Let $\widetilde{\omega} := \omega-d g_m$. Then $\widetilde{\omega}= \widetilde{f}_1 \, dx_1+\cdots+\widetilde{f}_{m-1} \, dx_{m-1}$, where for each $i$,
\[
\tilde{f}_i := f_i - \partial_i(g_m) = f_i -\partial_i(a_m) - \sum_{j=1}^{n} c_{m, j} \frac{\partial_i(b_{m, j})}{b_{m, j}} \in \overline{F}.
\]
This shows that the coefficient field is extended only algebraically. And the constant field is also extended algebraically, i.e., $\Const_{\partial_1,\dots,\partial_m}(\overline{F})  = \overline{C}$. Since $\omega$ is closed, so is $\widetilde{\omega}$, which implies that $\partial_m(\tilde f_i)=0$ for $i\in\{1, \ldots, m-1\}$. Thus $\widetilde{\omega}$ is a differential form that is free of $dx_m$ and has coefficients in $C_m(c_{m,1},\ldots,c_{m,n})$.
	
Repeating this process recursively shows that there are $a \in \overline{F}$, $c_j \in \overline{C}$ and nonzero $b_j\in \overline{C}F$ such that
\[
f_i = \partial_i(a) + \sum_{j} c_j \frac{\partial_i (b_j)}{b_j}
\]
for all $i \in \{1,\dots,m\}$.  This almost has the required form. The problem is that $a$ need not lie in $F$ and the $b_j$'s need not lie in $F(\{c_j\})$.  We have more work to do.
	
Let $A$ be the minimal Galois extension of $F$ generated by $a$, the $c_{j}$'s, and the $b_{j}$'s with $d := [A : F] < \infty$.   The trace map  $\operatorname{Tr}_{A/F} : A \to F$ is defined by 
$\operatorname{Tr}_{A/F}(a) = \sum_{\sigma \in G} \sigma(a)$ for $a\in A$, where $G$ be the Galois group of $A$ over $F$.   According to Theorem 3.2.4 in~\cite{Bronstein}, each $\sigma\in G$ commutes with $\partial_i$ for all $i  \in \{1,\ldots,m\}$.

Using the trace map, we obtain
\[
d \cdot f_i=\operatorname{Tr}_{A/F}(f_i) = \partial_i\bigl( \operatorname{Tr}_{A/F}(a)\bigr) + \sum_j \sum_{\sigma\in G} \sigma(c_j) \frac{\partial_i(\sigma(b_j))}{\sigma(b_j)}.
\]
Thus
\[
f_i = \partial_i\Bigl(\frac{\operatorname{Tr}_{A/F}(a)}{d}\Bigr) + \sum_j \sum_{\sigma\in G} \dfrac{\sigma(c_j)}{d} \frac{\partial_i(\sigma(b_j))}{\sigma(b_j)},	
\]
where
\[
\frac{\operatorname{Tr}_{A/F}(a)}{d} \in F, \quad \dfrac{\sigma(c_j)}{d} \in \overline{C} \quad \text{and} \quad  \sigma(b_j) \in \overline{C}F.
\]
This is closer to what we want, but the  $\sigma(b_j)$'s need not lie in $F(\{\sigma(c_j)\}_{j,\sigma})$.  So there is still more work to do.
	
Now let $L$ be the field extension of $F$ generated by $\sigma(c_j)$'s.  Note that $F \subseteq L \subseteq A$.  For the extension $L \subseteq A$ of degree $\ell:=[A:L]$, let $\operatorname{Tr}_{A/L} \colon A \rightarrow L$ and $N_{A/L} \colon A \rightarrow L$ be the trace and norm maps.  Recall that for any nonzero $b \in A$, according to Theorem~3.2.4 in~\cite{Bronstein}, we have
\[
 \operatorname{Tr}_{A/L}\Big(\frac{\partial_i(b)}{b}\Big) = \frac{\partial_i(N_{A/L}(b))}{N_{A/L}(b)}.
 \]
 Then applying  $\operatorname{Tr}_{A/L}$ to the above expression for $f_i$ yields
\begin{align*}
\ell \cdot f_i &= \operatorname{Tr}_{A/L} (f_i) = \ell \cdot \partial_i\Big(\frac{\operatorname{Tr}_{A/F}(a)}{d}\Big) + \sum_j \sum_{\sigma\in G} \dfrac{\sigma(c_j)}{d} \operatorname{Tr}_{A/L} \Big(\frac{\partial_i(\sigma(b_i))}{\sigma(b_i)}\Big)\\
&= \ell \cdot \partial_i\Big(\frac{\operatorname{Tr}_{A/F}(a)}{d}\Big) + \sum_j \sum_{\sigma\in G} \dfrac{\sigma(c_j)}{d} \frac{\partial_i(N_{A/L}(\sigma(b_i)))}{N_{A/L}(\sigma(b_i))},
\end{align*}
where the second line uses the above identity.  Hence
\[
f_i = \partial_i\Big(\frac{\operatorname{Tr}_{A/F}(a)}{d}\Big) + \sum_j \sum_{\sigma\in G} \dfrac{\sigma(c_j)}{d\cdot \ell} \frac{\partial_i(N_{A/L}(\sigma(b_i)))}{N_{A/L}(\sigma(b_i))},
\]
where $\operatorname{Tr}_{A/F}(a)/d \in F$, $ \dfrac{\sigma(c_j)}{d\cdot \ell} \in \overline{C}$, and $N_{A/L}(\sigma(b_j)) \in L = F(\{\sigma(c_j)\}_{j,\sigma})$. Since this holds for $i \in \{1,\dots,m\}$, it follows that $\omega = f_1\,dx_1 + \cdots +  f_m\,dx_m  = dg$, where 
\[
g = \frac{\operatorname{Tr}_{A/F}(a)}{d}  +  \sum_j \sum_{\sigma\in G} \dfrac{\sigma(c_j)}{d\cdot \ell} \log N_{A/L}(\sigma(b_i))
\]
has the required form. This completes the proof.
\end{proof}

The following example indicates the difficulty of extending Theorem~\ref{THM:Liouville1form} to $p$-forms for $p> 1$.

\begin{example}
Let $\omega = \bigl(\exp(x+y^2) + \exp(y+x^2)\bigr) \, dxdy$. It is easy to check that $\omega = d\eta$ with $\eta :=\exp(x+y^2) \, dy - \exp(y+x^2)\,dx$.  
		
	The recursive strategy used to prove Theorem~\ref{THM:1form} for 1-forms adapts nicely to prove Theorem~\ref{THM:pform} for $p$-forms.  However, this strategy fails when we try to adapt the proof of Theorem~\ref{THM:Liouville1form} to the above 2-form, which we write as  $\omega = f\, dxdy$ for $f := \exp(x+y^2) + \exp(y+x^2)$.  To see why, we follow the strategy of Example~\ref{EXAM:1/xy}.  The first step is to write $\omega = dy({-}f\,dx)$ and then find $g$ such that ${-}f = \partial_{y}(g)$.  But no elementary function has this property because $f$ contains $\exp(x+y^2)$.  Also, switching $x$ and $y$ does not help, since here $\omega = dx({}f\,dy)$, but there is no elementary function $g$ satisfying $f = \partial_{x}(g)$ because $f$ contains $\exp(y+x^2)$.
\end{example}

\section{Zeilberger's Creative Telescoping for Differential One-Forms with One Parameter}
\label{SECT:ZeilCT}
The method of \textit{creative telescoping} was first mentioned in van der Poorten's expository paper~\cite{Poorten1978} on Ap\'{e}ry's proof of the irrationality of $\zeta(3)$. Specifically, creative telescoping was used to show that the sum
\[ b_n := \sum_\ell b_{n,\ell} := \sum_{\ell=0}^n \binom{n}{\ell}^2\binom{n+\ell}{\ell}^2\]
satisfies the recurrence
\[(n+1)^3 b_{n+1} -(34n^3+51n^2+27n+5)b_{n} +n^3b_{n-1} =0, \text{ for all } n \geq 1.\]
The essential idea is to ``cleverly construct''~\cite{Poorten1978} the product
\[B_{n,\ell} =  4(2n+1)\big(\ell(2\ell+1)-(2n+1)^2\big) \cdot b_{n,\ell} \]
and observe that if $\Delta_\ell$ is the usual backward difference operator, then the identity
\begin{align*} 
(n+1)^3 b_{n+1,\ell} -(34n^3+51n^2+27n+5)b_{n,\ell} +n^3b_{n-1,\ell} &=B_{n,\ell}- B_{n,\ell-1}\\ &=: \Delta_\ell(B_{n,\ell}). 
\end{align*}
follows from a straightforward algebraic computation since $b_{n,\ell} =  \binom{n}{\ell}^2\binom{n+\ell}{\ell}^2$. Both $b_{n,\ell}$ and $B_{n,\ell}$ vanish automatically if $\ell < 0$ or $\ell > n$. For each fixed~$n \ge 1$, summing over all integers $\ell$ makes the right-hand side vanish, which proves the recurrence.

In general, for a given sequence $b_{n, \ell}$, the process of creative telescoping finds
a nontrivial linear recurrence $L\in k(n)\langle S_n\rangle$ and another operator $Q\in k(n, \ell)\langle S_n, S_{\ell}\rangle$ such that
\[
L(b_{n, \ell}) = \Delta_\ell(Q(b_{n, \ell})).
\]
In $k(n)\langle S_n\rangle$, $n$ is the operator of multiplication by $n$ and $S_n$ is the shift operator that maps $n$ to $n+1$.  The notation $k(n)\langle S_n\rangle$ reflects the fact that $n$ and $S_n$ do not commute since $S_n n = (n+1)S_n$. In the above identity, the operator $L$ is called a \textit{telescoper} for $b_{n, \ell}$ and $Q(b_{n, \ell})$ is a \textit{certificate} of $L$ for $b_{n, \ell}$. Zeilberger's fast algorithm~\cite{Zeilberger1990,Zeilberger1991} for creative telescoping was developed to solve this problem when restricted to so-called \textit{hypergeometric terms}  $b_{n ,\ell}$, meaning that both $S_n(b_{n,\ell})/b_{n,\ell}$ and $S_\ell(b_{n,\ell})/b_{n,\ell}$ are rational functions of $n$ and $\ell$.

In the differential analogue of creative telescoping~\cite{AlmkvistZeilberger1990},
sums and difference operators are replaced with integrals and differential operators, respectively.  Thus, given $F(x,y)$ in place of $b_{n,\ell}$, we seek
 a linear differential equation satisfied by the integral
\[f(x)= \int_{-\infty}^{+\infty} F(x,y) \, dy ,\]
where $F(x,y)$ is \textit{hyperexponential} with respect to $x,y$, i.e., $\partial_x(F(x,y))/F(x,y)$ and $\partial_y(F(x,y))/F(x,y)$ are both rational functions of $x$ and $y$.  In this case, there always exists a nonzero telescoper $L\in k(x)\left\langle \partial_x\right\rangle$ for $F(x,y)$ with respect to $y$ with a certificate~$G(x,y)$ such that
\[L(F(x,y))= \partial_y (G(x,y)).\]
(As above, the notation $k(x)\left\langle \partial_x\right\rangle$ reflects the fact that $x$ and $\partial_x$ do not commute since $\partial_x x = x\partial_x + 1$.) Then integrating with respect to $y$ leads to
\[L(f(x))=0.\]
In recent years, reduction-based creative telescoping~\cite{BostanChenChyzakLi, bostan-chen-chyzak-li-xin, BostanPierreBruno2013, ChenHuangKauersLi, ChenKauersKoutschan2016,bostan-chyzak-lairez-salvy, chen-vanhoeij-kauers-koutschan} has become a new trend because of its efficiency in computing telescopers without certificates. 

In this section, we assume that all integrals involve one parameter $t$ and the derivation with respect to $t$ is denoted by $\partial_t$. We will prove that nonzero telescopers exist for closed algebraic 1-forms and use Algorithm \texttt{HermiteOneForm} from Section~\ref{SECT:1form} to give a reduction-based creative telescoping algorithm to compute the telescoper for a closed rational $1$-form.

\begin{defin}
	Given a $1$-form $\omega = f_1 \hskip1pt dx_1+ \cdots + f_m \hskip1pt dx_m$ with $f_i$ algebraic over $k(t,x_1,\ldots,x_m)$, a linear differential operator $L$ in $k(t)\left<\partial_t\right>$ is called a \textit{telescoper} for $\omega$ if $L(\omega) := L(f_1) \, dx_1 + \cdots + L(f_m)\, dx_m =  dg$, where $g$ is algebraic  over $k(t,x_1,\ldots,x_m)$ and $d$ is the exterior derivative with respect to $x_1,\dots,x_m$. The function $g$ is called a \textit{certificate} of $L$.  A nonzero telescoper for $\omega$ of minimal degree in $\partial_t$ is called a \textit{minimal telescoper}.
\end{defin}

Observe that $L(\omega) = dg$ is equivalent to $L(f_i) = \partial_i(g)$ for each $i \in \{1,\dots,m\}$.  In the terminology of~\cite{ChenFengLiSinger2014}, 
$L$ is called a \textit{parallel telescoper}.

Creative telescoping for differential forms with D-finite coefficients has been studied in~\cite{ChenFengLiSinger2014, ChenFengLiSingerWatt2024}. There are two fundamental problems related to creative telescoping: one is the existence problem for telescopers and the other is the construction problem for telescopers. The following theorem concerning the existence of telescopers for differential 1-forms with algebraic-function coefficients is a special case of~\cite[Lemma~8]{ChenFengLiSinger2014} since algebraic functions are D-finite~\cite{Abel1827}.  Note that the proof of Lemma~8 in~\cite{ChenFengLiSinger2014} uses the same recursive argument as the proof of Theorem~\ref{THM:1form}. To be self-contained, we present the proof here in the language of differential forms. Moreover, the papers~\cite{ChenFengLiSinger2014, ChenFengLiSingerWatt2024} not only treat the existence problem of telescopers, but also present algorithms~\texttt{ParaTel} in~\cite{ChenFengLiSinger2014} to compute minimal telescopers for differential 1-form with hyperexponential-function coefficients and later~\texttt{Algorithm 14} in~\cite{ChenFengLiSingerWatt2024} for the case of general differential forms with D-finite coefficients, which are both based on the recursive idea used in the proof of Lemma~8. For a more special case of differential 1-form with rational-function coefficients, we first provide an alternative proof in Theorem~\ref{THM:CTrational} and then present a reduction-based algorithm~\texttt{CTOneForm} that differs from the recursive way in~\texttt{ParaTel}.

\begin{remark}
Even though the existence can also be derived from the finite-dimensionality of the de Rham cohomology groups of $D$-modules in the Bernstein class, endowed with the Gauss--Manin connection, the proofs of Theorems~\ref{THM:existCT1form} and \ref{THM:CTrational} given below have the same feature of being constructive 
as that in~\cite{ChenFengLiSinger2014, ChenFengLiSingerWatt2024}, so that we can design efficient algorithms for computing telescopers. 
\end{remark}

\begin{thm}\label{THM:existCT1form}
	Let $\omega= f_1 \, dx_1+ \cdots + f_m \, dx_m$ be a $1$-form with coefficients algebraic over $k(t,x_1,\ldots,x_m)$. If $\omega$ is closed, then it always has a nonzero telescoper.
\end{thm}

\begin{proof}
We proceed by induction on $m$, the number of variables $x_1,\ldots,x_m$. For $m=1$, one can construct a nonzero telescoper for $f_1$ with respect to $x_1$ by the method in~\cite{ChenKauersSinger2012,ChenKauersKoutschan2016}.
Suppose now that $m>1$ and that the theorem holds for the case when the number of variables is less than $m$. For $f_m \in \overline{k(t,x_1,\ldots,x_{m})}$, there exists a nonzero telescoper $L_m \in k(t,x_1,\ldots,x_{m-1})\left\langle \partial_t \right\rangle$ with respect to $x_m$ and a certificate $g_m$ such that
\[L_m(f_m)=\partial_m(g_m).\]

We claim that $L_m$ can be chosen to be free of $x_1,\ldots,x_{m-1}$.  Without loss of generality, assume that all coefficients of $L_m$ lie in the polynomial ring $k[t,x_1,\ldots,x_{m-1}]$, otherwise we can multiply {both sides of the above equation by the least common multiple of the denominators}. Among such $L_m$, we first show that $L_m$ can be chosen to be free of $x_1$.  If this is not the case, then 
pick $L_m$ with $\deg_{x_1}(L_m)$ minimal, denoted by $n$, and note that $n > 0$.  Then write
	\[L_m = x_1^{n}L_{m,n}+\cdots+x_1L_{m,1}+L_{m,0}, \]
where $L_{m,j} \in k[t,x_2,\ldots,x_{m-1}]\left\langle \partial_t\right\rangle$ for each $j \in \{1,\ldots,n\}$.  To make the following computation clearer, we use ``$\circ$'' to represent the composition of operators. Observe that
\[
\partial_1\circ x_1^j = x_1^j \circ \partial_1 + j x_1^{j-1}.
\]
Thus
\[
\partial_1\circ L_m = x_1^n \circ \partial_1 \circ L_{m,n} + \cdots + x_1 \circ \partial_1 \circ L_{m,1}+ \partial_1 \circ L_{m,0}+ \sum_{j=1}^{n} j x_1^{j-1}\circ L_{m,j}.
\]
Since $\omega$ is closed, we have $\partial_1(f_m)= \partial_m (f_1)$. Thus, for each $j \in \{1, \ldots, n\}$,
\[
x_1^j \circ \partial_1 \circ L_{m,j} (f_m) = x_1 ^{j}\circ L_{m,j} \circ \partial_1 (f_m)= x_1 ^{j}\circ L_{m,j} \circ \partial_m (f_1)= \partial_m (x_1^jL_{m,j}(f_1)).
\]
Now define $\tilde{L}_m := \partial_1 \circ L_{m,0}+ \sum_{j=1}^{n} j x_1^{j-1}\circ L_{m,j}$.  Putting this all together and using $\partial_1(L_m(f_m))= \partial_1(\partial_m(g_m)) = \partial_m(\partial_1(g_m))$, we obtain
\[
\tilde{L}_m (f_m)=\partial_m\Big(\partial_1(g_m)-\sum_{j=1}^{n}x_1^j L_{m,j}(f_1)\Big),
\]
which implies that $\tilde{L}_m$ is a telescoper for $f_m$ with respect to $x_m$ with the degree {in $x_1$ at most} $n-1$, leading to a contradiction. This shows that $L_m$ can be chosen to be free of $x_1$.   Then, among all $L_m$'s that are free of $x_1$, a similar argument using $x_2$ shows that $L_m$ can be chosen to be free of $x_1$ and $x_2$.  Continuing in this way, we see that 
our claim is true, i.e., $L_m \in k(t)\left\langle \partial_t\right\rangle$.
	
Since $\partial_i\circ L_m = L_m \circ \partial_i$ for each $i \in \{1,\ldots,m\}$, we see that $L_m(\omega)$ is closed with coefficients algebraic over $k(t,x_1,\ldots,x_m)$. Applying $L_m$ to $\omega$ leads to
\[
L_m (\omega)= L_m(f_1) \, dx_1 +\cdots +  L_m(f_{m-1})\, dx_{m-1} + \partial_m (g_m) \, dx_m.
\]
By subtracting $dg_m$, we see that $L_m(\omega) - dg_m$ is a closed algebraic differential form free of $x_m$ and~$dx_m$. By induction, there exists a nonzero telescoper $\tilde{L}$ with a certificate $\tilde{g}$ such that
\[
\tilde{L} (L_m(\omega) - dg_m) = d \tilde{g},
\]
which implies that
\[
\tilde{L}\circ L_m (\omega)= d\big(\tilde{L}(g_m)+\tilde{g}\big).
\]
Thus $L := \tilde{L}\circ L_m$ is a nonzero telescoper for $\omega$ with a certificate $g = \tilde{L}(g_m)+\tilde{g}$.
\end{proof}

Considering rational closed $1$-forms, we next present an even simpler proof of the existence of telescopers based on Theorem~\ref{THM:1form} and Abel's classical lemma~\cite{Abel1827,BCSLS07} that any algebraic function satisfies a linear differential equation with polynomial coefficients.

\begin{thm}\label{THM:CTrational}
    Let $\omega= f_1 \, dx_1+ \cdots + f_m \, dx_m$ be a $1$-form with~$f_i \in k(t,x_1,\ldots,x_m)$. If $\omega$ is closed, then it always has a nonzero telescoper.
\end{thm}
\begin{proof}
    From Theorem~\ref{THM:1form}, we have
    \[
    f_i = \partial_i(a) + \sum_{j} c_j \frac{\partial_i(b_j)}{b_j}, \ \text{for } i\in\{1,\ldots,m\},
    \]
    where~$a\in k(t,x_1,\ldots,x_m)$, $c_j\in \overline{k(t)}$, and~$b_j \in k(t)(c_j)(x_1,\ldots,x_m)$. Since $\partial_t$ commutes with~$\partial_i$, we have 
\begin{equation}
    \begin{aligned} \label{EQ:reduct1form}
        \partial_t(f_i) &= \partial_i\bigl(\partial_t(a)\bigr) + \sum_j \Biggl(\! c_j \partial_i \biggl(\frac{\partial_t(b_j)}{b_j} \biggr)+ 
        \partial_t(c_j)\frac{\partial_i(b_j)}{b_j}\Biggr)\\
        &= \partial_i \biggl(\!\partial_t(a) + \sum_j c_j \frac{\partial_t(b_j)}{b_j} \biggr) + \sum_j \partial_t(c_j)\frac{\partial_i(b_j)}{b_j}.
    \end{aligned}
    \end{equation}
    By Abel's lemma, for each $j$, there exists a nonzero $L_j \in k(t)\left\langle \partial_t \right\rangle$ such that $L_j(c_j)=0$.
    Using standard algorithms, compute the least common left multiple $L$ of the $L_j$'s.  Then $L$ is a nonzero operator  such that $L(c_j)=0$ for each~$j$.
    Hence, applying \eqref{EQ:reduct1form} iteratively leads to 
    \[
        L(f_i) = \partial_i(g) + \sum_j L(c_j) \frac{\partial_i(b_j)}{b_j}= \partial_i(g), \text{ where }g\in k(t)(\{c_j\})(x_1,\ldots,x_m).
    \]
    Using the trace map similarly as in the proof of Theorem~\ref{THM:Liouville1form}, we have 
    \[
        L(f_i) = \partial_i(g) , \text{ for some }g \in k(t,x_1,\ldots,x_m). 
    \]
    Since the above equation holds for every~$i$, the proof is complete.
\end{proof}

We now present a reduction-based algorithm for computing minimal telescopers for 
 closed rational 1-forms which generalizes the method in~\cite{BostanChenChyzakLi}
 from rational functions to rational 1-forms. The correctness of this algorithm follows from the above existence theorem on telescopers for closed rational 1-forms
 and the properties of Hermite reduction for these forms in Section~\ref{SECT:1form}.

\medskip
\noindent
	\textbf{Algorithm} \texttt{CTOneForm:}
	Given a closed rational $1$-form $\omega = f_1 \, dx_1+\cdots+ f_m \, dx_m$, where $f_i \in k(t,x_1,\ldots,x_m)$, return a minimal telescoper $L\in k(t)\left\langle \partial_t\right\rangle$ and a certificate $g$ of $L$.
\begin{enumerate}
\item Decompose $\omega$ as $\omega= dg_0+ r_0$ by \texttt{HermiteOneForm} from Section~\ref{SECT:1form}.
\item Set $n:=0$ and~$\mathcal{L} := [r_n]$.
\item 
While the elements of~$\mathcal{L}$ are linearly independent over~$k(t)$ do \\
\indent Set~$n:=n+1$. Apply \texttt{HermiteOneForm} to $\partial_t(r_{n-1})$:
\[
    \partial_t(r_{n-1}) = dg_{n} + r_{n}.
\]
\indent Append $r_{n}$ to~$\mathcal{L}$.
\item Find a nontrivial linear relation~$\sum_{i=0}^n a_i(t) r_i = 0$ where~$a_i(t) \in k(t)$.
\item Return~$L = \sum_{i=0}^n a_i(t) \partial_t^i$ and~$g =\sum_{i=0}^n a_i(t) \sum_{j=0}^i \partial_t^{i-j}\bigl(g_{j}\bigr)$.
\end{enumerate}

\begin{example}\label{EXAM:Telescoper}
    Let $\omega$ be as the same $1$-form as in Example~\ref{EXAM:Hermite}. Repeatedly applying \texttt{HermiteOneForm} gives
    \begin{gather*}
    g_0 = \frac{1}{xyz},\quad r_0:= \frac{t}{x} \, dx+ \frac{t}{y} \, dy + \frac{t^2+1}{z}\, dz,\\
    g_1 = 0,\quad r_1:= \frac{1}{x} \, dx+ \frac{1}{y} \, dy + \frac{2t}{z}\, dz,\\
    g_2 = 0,\quad r_2:= 0 \, dx+ 0 \, dy + \frac{2}{z}\, dz.
    \end{gather*}
    Note that $\{ r_0, r_1\}$ is linearly independent over~$k(t)$, whereas $\{ r_0, r_1, r_2\}$ is linearly dependent, as shown by
    \[
        (t^2-1) r_2 -2tr_1 + 2 r_0 =0.
    \]
    So a minimal telescoper for~$\omega$ is 
    \[
        L = (t^2-1)\partial_t^2 - 2t \partial_t + 2 ,
    \]
    with a certificate  
    \[
        g := (t^2-1)\cdot \partial_t^2(r_0) -2t \cdot \partial_t(r_0) + 2 \cdot r_0 = \frac{2}{xyz}.
    \]
\end{example}

\section*{Acknowledgements}

The authors are very grateful to Michael Singer for making this joint work possible and Pierre Lairez for many discussions about Picard's problem. 
We would also like to thank the anonymous referees for their comments, which helped us significantly improve the paper.

\section*{Declaration on the Use of Generative AI}

During the preparation of this manuscript, the authors used ChatGPT to assist with language polishing. The generated content was critically reviewed, edited, and verified by the authors. The authors take full responsibility for the accuracy, originality, and integrity of the final manuscript.

\bibliographystyle{abbrv}
\bibliography{diffform}
\end{document}